\documentclass[MSNbibl,number,citesort,seceqn,dvips]{arxbj}

\aid{0}
\volume{20}
\issue{2}
\pubyear{2014}
\firstpage{486}
\lastpage{513}
\doi{10.3150/12-BEJ494} 

\makeatletter
\newcommand{\RMo}{\mathrm{o}}

\newcommand{\RMe}{\mathrm{e}}

\newcommand{\mrmd}{\mathrm{d}}

\newcommand{\iint}{\int\!\!\int}

\newcommand{\I}{{\mathbf1}}

\newtheorem{proposition}{Proposition}[section]
\newtheorem{theorem}[proposition]{Theorem}
\newtheorem{corollary}[proposition]{Corollary}
\newtheorem{lemma}[proposition]{Lemma}

\newremark{remark}[proposition]{Remark}

\newremark{example}[proposition]{Example}

\newcommand{\R}{{\mathbb R}}
\newcommand{\bS}{{\mathbb S}^{d-1}}
\newcommand{\N}{{\mathbb N}}
\newcommand{\BP}{\mathbb{P}}
\newcommand{\BV}{\mathbb{V}}
\newcommand{\BE}{\mathbb{E}}
\newcommand{\BQ}{\mathbb{Q}}
\newcommand{\bD}{{\mathbf D}}
\newcommand{\bF}{{\mathbf F}}
\newcommand{\bN}{{\mathbf N}}
\newcommand{\BX}{{\mathbb X}}
\newcommand{\BY}{{\mathbb Y}}
\newcommand{\bg}{{\mathbf g}}

\makeatother

\begin{document}
\begin{frontmatter}

\title{Perturbation analysis of Poisson processes}
\runtitle{Perturbation analysis of Poisson processes}

\begin{aug}
\author{\fnms{G\"unter} \snm{Last}\corref{}\ead[label=e1]{guenter.last@kit.edu}}
\runauthor{G. Last} 
\address{Karlsruhe Institute of Technology,
Institut f\"ur Stochastik, Kaiserstra{\ss}e 89, D-76128 Karlsruhe,
Germany. \printead{e1}}
\end{aug}

\received{\smonth{3} \syear{2012}}
\revised{\smonth{10} \syear{2012}}

%
\begin{abstract}
We consider a Poisson process $\Phi$ on a general phase space. The
expectation of a function of $\Phi$ can be considered as a functional
of the intensity measure $\lambda$ of $\Phi$. Extending earlier results
of Molchanov and Zuyev [\textit{Math. Oper. Res.} \textbf{25} (2010)
485--508] on finite Poisson processes, we study the
behaviour of this functional under signed (possibly infinite)
perturbations of $\lambda$. In particular, we obtain general
Margulis--Russo type formulas for the derivative with respect to
non-linear transformations of the intensity measure depending on some
parameter. As an application, we study the behaviour of expectations of
functions of multivariate L\'evy processes under
perturbations of the L\'evy measure. A key ingredient of our approach
is the explicit Fock space representation obtained in Last and Penrose
[\textit{Probab. Theory Related Fields} \textbf{150} (2011)
663--690].
\end{abstract}

%

\end{frontmatter}

\section{Introduction}\label{secintro}

The aim of this paper is to advance
the perturbation analysis of a \textit{Poisson process} $\Phi$
on a general measurable space $(\BX,\mathcal{X})$.
For any $\sigma$-finite measure $\lambda$ on $(\BX,\mathcal{X})$,
we let $\Pi_\lambda$ denote the distribution of a Poisson
process with \textit{intensity measure} $\lambda$, see, for example,
\cite{Kallenberg}, Chapter 12.
Further we let $\BP_\lambda$ be a probability measure on some fixed
measurable sample space such that $\BP_\lambda(\Phi\in\cdot)=\Pi_\lambda$.
We let $\BE_\lambda$ denote the expectation operator with respect
to $\BP_\lambda$.
Let $f(\Phi)$ be some (measurable) function of $\Phi$.
Under certain assumptions on $f$,
Molchanov and Zuyev \cite{MolZu00}
showed for finite measures $\lambda$ and $\nu$ the
variational formula
%
\begin{equation}
\label{main1} \BE_{\nu} f(\Phi)=\BE_\lambda f(\Phi) +\sum
^\infty_{n=1}\frac{1}{n!}\int\bigl(
\BE_\lambda D^n_{x_1,\ldots,x_n} f(\Phi)\bigr) (\nu-
\lambda)^n\bigl(\mrmd(x_1,\ldots,x_n)\bigr),
\end{equation}
where
%
\begin{equation}
\label{Dsymmetric0} D^n_{x_1,\ldots,x_n}f(\Phi)= \sum
_{J
\subset\{1,2,\ldots,n\}}(-1)^{n-|J|} f \biggl(\Phi+\sum
_{j\in
J}\delta_{x_j} \biggr),\qquad x_1,\ldots,x_n\in\BX, n\in\N.
\end{equation}
Here, $|J|$ denotes the number of elements of $J$,
while $\delta_{x}$ is the Dirac measure located at a point $x\in\BX$.
It is common to say that $\nu$ results from $\lambda$
by adding the \textit{perturbation} $\nu-\lambda$.

In this paper, we shall extend (\ref{main1})
to $\sigma$-finite measures $\lambda$ and $\nu$.
One can use a pathwise defined thinning and superposition
construction to move from $\BP_\lambda$ to $\BP_\nu$, see Remark
\ref{thsup}.
In general, $\nu-\lambda$ is a \textit{signed measure} that cannot be defined
on the whole $\sigma$-field $\mathcal{X}$. Integration with respect to
$\nu-\lambda$, however, is well defined via (\ref{signed}) below.
Under an additional assumption on $\nu$ and $\lambda$
(satisfied for positive, negative and many other perturbations of
$\lambda$),
we shall establish a condition that
is necessary and sufficient for (\ref{main1}) to
hold for all bounded functions of $\Phi$.
If, for instance, $\lambda\le\nu$, this
condition is equivalent to the absolute continuity
$\Pi_\lambda\ll\Pi_\nu$. The variational formula does not only
hold for bounded
functions but under a more general second moment assumption on $f$.

A consequence of (\ref{main1}) are derivative formulas of the form
%
\begin{equation}
\label{deriv00} \frac{\mrmd}{\mrmd \theta}\BE_{\lambda+\theta(\nu-\lambda)} f(\Phi) =\int
\BE_{\lambda+\theta(\nu-\lambda)} D_xf(\Phi) (\nu-\lambda ) (\mrmd x), \qquad\theta\in[0,1],
\end{equation}
where $D_x:=D^1_x$ is the first order difference (or add one cost)
operator. This can be generalized to non-linear perturbations of
$\lambda$
and to more than one parameter.
Such formulas are useful in the performance evaluation, optimization
and simulation of discrete
event systems \cite{HoCao91,AsGlynn07}. Applications in a
spatial setting can be found in \cite{BaKlZu95,MolZu00}.
Equation (\ref{deriv00}) can be seen as a
Poisson version of the \textit{Margulis--Russo formula}
for Bernoulli random fields (see, e.g., \cite{BoRiordan06}).
Such formulae are, for instance, an important tool in both
discrete and continuum percolation theory.

The extension of the identity (\ref{main1}) from finite
to $\sigma$-finite measures is a non-trivial task.
Our approach is based on a combination
of the recent Fock space representation in \cite{LaPe11} with
classical results in \cite{Brown71} on the absolute continuity
of Poisson process distributions. A related approach to
derivatives of the type (\ref{deriv00}) for marked point processes on
the real line
was taken in \cite{Decreuse98}.
For Poisson processes on the line
and under a (rather strong) continuity assumptions on $f$ the result
(\ref{main1}) can be considered as a special case
of the main result in \cite{Bl95}.

The paper is organized as follows. In Section \ref{secpre}, we introduce
some basic notation and recall facts about the
Fock space representation and likelihood functions
of Poisson processes. In Section \ref{secfinite}, we use an
elementary but illustrative argument to prove a simple version of
(\ref{main1}). In Section \ref{secinf}, we prove and discuss
Theorem \ref{tvar3}, which is the main result of this paper.
In Section \ref{secnec}, we derive conditions on $\lambda$ and
$\nu$ that are necessary for (\ref{main1}) to hold for all
bounded functions $f$. In some cases these conditions
are also sufficient. Section \ref{secder} gives general
Margulis--Russo type formulas for derivatives.
The final Section \ref{seclevy} treats perturbations
of the L\'evy measure of a L\'evy process in $\R^d$.

\section{Preliminaries}\label{secpre}

Let $\bN$ be the space of integer-valued $\sigma$-finite measures
$\varphi$
on $\BX$ equipped with the smallest $\sigma$-field $\mathcal{N}$
making the mappings $\varphi\mapsto\varphi(B)$ measurable for all
$B\in\mathcal{X}$. We fix a measurable
mapping $\Phi\dvtx \Omega\rightarrow\bN$, where
$(\Omega,\mathcal{F})$ is some abstract measurable (sample) space.
For any $\sigma$-finite measure $\lambda$ on $(\BX,\mathcal{X})$
we let $\BP_\lambda$ be a probability measure on
$(\Omega,\mathcal{F})$ such that $\BP_\lambda(\Phi\in\cdot)=\Pi_\lambda$
is the distribution of a Poisson process with intensity measure
$\lambda$.

For any measurable $f\dvtx \bN\rightarrow\R$ and $x\in\BX$ the function
$D_xf$ on $\bN$ is defined by
%
\begin{equation}
\label{addone} D_{x}f(\varphi):=f(\varphi+\delta_{x})-f(
\varphi), \qquad\varphi\in \bN.
\end{equation}
The \textit{difference operator} $D_x$ and its iterations
play a central role in the variational analysis of Poisson processes.
For $n\ge2$ and $(x_1,\ldots,x_n)\in\BX^n$
we define a function
$D^{n}_{x_1,\ldots,x_n}f\dvtx \bN\rightarrow\R$ inductively by
%
\begin{equation}
\label{differn} D^{n}_{x_1,\ldots,x_{n}}f:=D^1_{x_{1}}D^{n-1}_{x_2,\ldots,x_{n}}f,
\end{equation}
where $D^1:=D$ and $D^0f = f$.
Note that
%
\begin{equation}
\label{Dsymmetric} D^n_{x_1,\ldots,x_n}f(\varphi)= \sum
_{J \subset\{1,2,\ldots,n\}}(-1)^{n-|J|} f \biggl(\varphi+\sum
_{j\in J}\delta_{x_j} \biggr),
\end{equation}
where $|J|$ denotes the number of elements of $J$. This
shows that the operator $D^n_{x_1,\ldots,x_n}$
is symmetric in $x_1,\ldots,x_n$, and that
$(\varphi,x_1,\ldots,x_n)\mapsto D^n_{x_1,\ldots,x_n}f(\varphi)$
is measurable.

From \cite{LaPe11}, Theorem 1.1, we obtain for any
measurable $f,g\dvtx \bN\rightarrow\R$ satisfying $\BE_\lambda f(\Phi
)^2<\infty$
and $\BE_\lambda g(\Phi)^2<\infty$ that
%
\begin{equation}
\label{28} \BE_\lambda f(\Phi)g(\Phi)= \sum^\infty_{n=0}
\frac{1}{n!}\int\bigl(\BE_\lambda D^n_{x_1,\ldots,x_n} f(
\Phi)\bigr) \bigl(\BE_\lambda D^n_{x_1,\ldots,x_n} g(\Phi)\bigr)
\lambda^n\bigl(\mrmd(x_1,\ldots,x_n)\bigr),
\end{equation}
where the summand for $n=0$ has to be interpreted as
$(\BE_\lambda f(\Phi))(\BE_\lambda g(\Phi))$. (The integral of a constant
$c$ with respect to $\lambda^0$ is interpreted as $c$.)

Next, we recall a result from \cite{Brown71} in a slightly
modified form. Consider two $\sigma$-finite measures
$\nu,\rho$ on $\BX$ such that $\nu\ll\rho$, that is,
$\nu$ is absolutely
continuous with respect to $\rho$. Let
$h:=\mrmd \nu/\mrmd \rho$
be the corresponding density (Radon--Nikodym derivative) and assume that
%
\begin{equation}
\label{hsquare} \int(h-1)^2\,\mrmd \rho<\infty.
\end{equation}
This implies that the sets $C_n:=\{|h-1|\ge1/n\}$, $n\in\N$,
have finite measure with respect to both $\nu$ and $\rho$,
cf. also \cite{Brown71}.
Define measurable functions $L_n\dvtx \bN\rightarrow[0,\infty)$
by
%
\begin{equation}
\label{Ln} L_n(\varphi):=\I\bigl\{\varphi(C_n)<\infty
\bigr\} \RMe ^{\rho(C_n)-\nu(C_n)}\prod_{y\in\varphi_{C_n}}h(y),
\end{equation}
where $\varphi_B$ is the restriction of $\varphi\in\bN$ to a measurable
set $B\subset\BX$ and
the product is over all points of the support
of $\varphi_{C_n}$ taking into account the multiplicities, that is,
\[
\prod_{y\in\varphi_{C_n}}h(y):= \exp \biggl[\int
_{C_n}\ln h(y)\varphi(\mrmd y) \biggr],
\]
where $\ln0:=-\infty$.\vadjust{\goodbreak}
%
\begin{proposition}\label{pbrown} With $\nu$ and $\rho$ as above
we have for any measurable $g\dvtx \bN\rightarrow\R$ that
%
\begin{equation}
\label{ac} \BE_{\nu}g(\Phi) =\BE_{\rho}L_{\nu,\rho}(\Phi)g(
\Phi),
\end{equation}
where
%
\begin{equation}
\label{L} L(\varphi):=L_{\nu,\rho}(\varphi):=\liminf_{n\to\infty
}L_n(
\varphi)
\end{equation}
if this limit inferior is finite and $L(\varphi):=L_{\nu,\rho
}(\varphi):=0$
otherwise. Furthermore,
%
\begin{equation}
\label{Lsquare} \BE_{\rho}L_{\nu,\rho}(\Phi)^2<\infty.
\end{equation}
\end{proposition}
\begin{pf}
It follows as in the proof of Theorem 1 in
\cite{Brown71} that $L_n(\Phi)$ converges
$\BP_{\rho}$-a.s. to a random variable $Y$
such that $\BE_{\nu}g(\Phi) =\BE_{\rho}Y g(\Phi)$
for all measurable $g\dvtx \bN\rightarrow\R$.
Hence, (\ref{ac}) holds. Furthermore, we have for any
$n\in\N$ that
\begin{eqnarray*}
\BE_{\rho}L_n(\Phi)^2&=& \exp\bigl[2
\rho(C_n)-2\nu(C_n)\bigr]\BE_{\rho} \biggl[\prod
_{y\in\Phi\cap
C_n}h(y)^2 \biggr]
\\
&=&\exp \biggl[2\rho(C_n)-2\int_{C_n}h \,\mrmd \rho
\biggr]\exp \biggl[\int_{C_n}\bigl(h^2-1\bigr)\,\mrmd \rho
\biggr]
\\
&=&\exp \biggl[\int_{C_n}(h-1)^2\,\mrmd \rho \biggr],
\end{eqnarray*}
where we have used a well-known property of Poisson processes
to obtain the second equality.
(Because $\rho(C_n)<\infty$ one can use a direct calculation
based on the mixed sample representation
or take $f:=-\ln h^2$ in \cite{Kallenberg}, Lemma 12.2(i),
see also \cite{MKM82}, 1.5.6.)
Fatou's lemma implies that
\[
\BE_{\rho}L(\Phi)^2\le\exp \biggl[\int(h-1)^2\,\mrmd
\rho \biggr],
\]
which is finite by assumption (\ref{hsquare}).
\end{pf}
%
\begin{remark}\label{remarkac} As noted above,
(\ref{hsquare}) implies that $\Pi_\nu\ll\Pi_{\rho}$.
The converse is generally not true.
However, if $h$ is bounded then (\ref{hsquare}) is necessary
and sufficient for $\Pi_\nu\ll\Pi_{\rho}$.
This follows from the main result in \cite{Brown71},
see also Theorem 1.5.12 in \cite{MKM82}.
\end{remark}

\section{Finite non-negative perturbations}\label{secfinite}

In this section, we fix a $\sigma$-finite measure $\lambda$
on $\BX$ and a finite measure $\mu$ on $\BX$.
In this case, we can derive the variational formula (\ref{main1})
for $\nu:=\lambda+\mu$ under a minimal integrability assumption on the
function $f$. Our proof (basically taken from \cite{MolZu00})
is elementary but instructive.
%
\begin{theorem}\label{tvar} Let
$f\dvtx \bN\rightarrow\R$ be a measurable function such that
$\BE_{\lambda+\mu}|f(\Phi)|<\infty$. Then (\ref{main1}) holds,
where all expectations exist and the series converges absolutely.
\end{theorem}
\begin{pf}
We perform a formal calculation using Fubini's theorem.
This will be justified below.
Denoting the right-hand side of (\ref{main1}) by $I$, and
using (\ref{Dsymmetric0}), we have that
\[
I=\sum^\infty_{n=0}\frac{1}{n!} \int
\biggl(\sum_{J \subset\{1,\ldots,n\}}(-1)^{n-|J|}
\BE_\lambda f\biggl(\Phi+\sum_{j\in J}
\delta_{x_j}\biggr) \biggr)\mu^n\bigl(\mrmd(x_1,\ldots,x_n)\bigr).
\]
By symmetry,
\begin{eqnarray*}
I&=&\sum^\infty_{n=0}\frac{1}{n!}\sum
^{n}_{m=0}\pmatrix{n
\cr
m}(-1)^{n-m}
\mu(\BX)^{n-m} \int\BE_\lambda f(\Phi+\delta_{x_1}+
\cdots+\delta_{x_m}) \mu^m\bigl(\mrmd(x_1,\ldots,x_m)\bigr)
\\
&=&\sum^{\infty}_{m=0}\frac{1}{m!}\sum
^\infty_{n=m}\frac
{(-1)^{n-m}}{(n-m)!} \mu(
\BX)^{n-m} \int\BE_\lambda f(\Phi+\delta_{x_1}+\cdots+
\delta_{x_m}) \mu^m\bigl(\mrmd(x_1,\ldots,x_m)\bigr)
\\
&=&\RMe^{-\mu(\BX)}\sum^{\infty}_{m=0}
\frac{1}{m!} \int\BE_\lambda f(\Phi+\delta_{x_1}+\cdots+
\delta_{x_m}) \mu^m\bigl(\mrmd(x_1,\ldots,x_m)\bigr)
\\
&=&\int\BE_\lambda f(\Phi+\varphi)\Pi_\mu(\mrmd \varphi),
\end{eqnarray*}
where in the last step we have used the mixed sample representation
of finite Poisson processes, see, for example, \cite{Kallenberg}, Theorem 12.7.
Noting that the distribution $\BP_{\lambda+\mu}(\Phi\in\cdot)$ is that
of a sum of two independent Poisson processes with intensity measures
$\lambda$ and $\mu$, respectively, we obtain that $I=\BE_{\lambda
+\mu} f(\Phi)$,
as desired.

To justify the use of Fubini's theorem, we need to show that
\[
c:=\sum^\infty_{n=0}\frac{1}{n!}\sum
^{n}_{m=0}\pmatrix{n
\cr
m}\mu (
\BX)^{n-m} \int\BE_\lambda\bigl|f(\Phi+\delta_{x_1}+\cdots+
\delta_{x_m})\bigr| \mu^m\bigl(\mrmd(x_1,\ldots,x_m)\bigr)<\infty.
\]
By a similar calculation as above,
\begin{eqnarray*}
c&=&\RMe^{\mu(\BX)}\sum^{\infty}_{m=0}
\frac{1}{m!} \int\BE_\lambda\bigl|f(\Phi+\delta_{x_1}+\cdots+
\delta_{x_m})\bigr| \mu^m\bigl(\mrmd(x_1,\ldots,x_m)\bigr)
\\
&=&\RMe^{2\mu(\BX)}\BE_{\lambda+\mu} \bigl|f(\Phi)\bigr|<\infty.
\end{eqnarray*}
This proves the theorem.
\end{pf}

\section{General perturbations}\label{secinf}

In this section, we allow also signed and infinite perturbations
of the intensity measure of $\Phi$. This requires
more advanced techniques, as the Fock space
representation (\ref{28}) and Proposition \ref{pbrown}.

We consider two $\sigma$-finite measures $\lambda$ and $\nu$ on $\BX$.
We take a $\sigma$-finite measure $\rho$ dominating
$\lambda$ and $\nu$, that is, $\lambda+\nu\ll\rho$.
Let $h_\lambda:=\mrmd \lambda/\mrmd \rho$, $h_\nu:=\mrmd \nu/\mrmd \rho$.
The integral of a measurable function
$g\dvtx \BX^n\rightarrow\R$ with respect to $(\nu-\lambda)^n$
is defined by
%
\begin{equation}
\label{signed} \int g\,\mrmd(\nu-\lambda)^n:=\int g(x_1,\ldots,x_n)
(h_\nu-h_\lambda)^{\otimes
n}(x_1,\ldots,x_n)\rho^n\bigl(\mrmd(x_1,\ldots,x_n)\bigr),
\end{equation}
where, for any function $h\dvtx \BX\rightarrow\R$, the function
$h^{\otimes n}\dvtx \BX^n\rightarrow\R$ is given by
\[
h^{\otimes n}(x_1,\ldots,x_n):=\prod
^n_{j=1}h(x_j).
\]
Note that our definition of $\int g\,\mrmd(\nu-\lambda)^n$
does not depend on the choice of $\rho$.
The following theorem is the main result of this paper.
%
\begin{theorem}\label{tvar3}
Assume that
%
\begin{equation}
\label{22a} \int(1-h_\lambda)^2 \,\mrmd \rho+ \int(1-h_\nu)^2
\,\mrmd \rho<\infty.
\end{equation}
Let $f\dvtx \bN\rightarrow\R$ be a measurable function such that
$\BE_{\rho}f(\Phi)^2<\infty$.
Then
%
\begin{equation}
\label{main2} \sum^\infty_{n=0}
\frac{1}{n!}\int\bigl|\BE_\lambda D^n_{x_1,\ldots,x_n}f(\Phi)\bigr|
|h_\nu-h_\lambda|^{\otimes n}(x_1,\ldots,x_n)\rho^n\bigl(\mrmd(x_1,\ldots,x_n)\bigr)<\infty
\end{equation}
and (\ref{main1}) holds.
\end{theorem}
%
\begin{remark}\label{thsup} Let $\Phi_\lambda$ be a Poisson process
with intensity measure $\lambda$ defined on some
abstract probability space. Then we can use independent
thinning and superposition to generate
a Poisson process $\Phi_\nu$ with intensity measure $\nu$.
Let $A:=\{x\in\BX\dvtx h_\lambda(x)>h_\nu(x)\}$
and define $p\dvtx \BX\rightarrow[0,1]$ by $p(x):=h_\nu(x)/h_\lambda(x)$
for $x\in A$ and by $p(x):=1$, otherwise.
Let $\Phi'$ be a $p$-thinning of $\Phi_\lambda$, see
\cite{Kallenberg}, Chapter 12. Then $\Phi'$ is a Poisson
process with intensity measure
\[
p(x)\lambda(\mrmd x)=\I_A(x)h_\nu(x)\rho(\mrmd x)+
\I_{\BX\setminus
A}h_\lambda(x)\rho(\mrmd x).
\]
Let $\Phi''$ be a Poisson process with
intensity measure $\I_{\BX\setminus A}(x)(h_\nu(x)-h_\lambda
(x))\rho(\mrmd x)$,
independent of $\Phi'$.
Then $\Phi'+\Phi''$ is a Poisson process with
intensity measure $h_\nu(x)\rho(\mrmd x)=\nu(\mrmd x)$.
In some applications, it might be convenient to couple $\Phi_\lambda$
and the perturbed process $\Phi_\nu$ in a different way.
For instance, $\Phi_\lambda$ could be an independent marking
of a homogeneous Poisson process of arrival times and
  one might wish to keep the times and to change only the marks.
\end{remark}

\begin{pf*}{Proof of Theorem \ref{tvar3}} By assumption (\ref{22a}), we can apply
Proposition \ref{pbrown} to both $\lambda$ and $\nu$.
It follows from (\ref{Ln}) and (\ref{L}) that
$L_{\lambda,\rho}(\Phi+\delta_x)=h_\lambda(x)L_{\lambda,\rho
}(\Phi)$
for all $x\in\BX$. Therefore,
\[
D^n_{x_1,\ldots,x_n} L_{\lambda,\rho}(\Phi)= L_{\lambda,\rho}(\Phi)
\prod^n_{i=1}\bigl(h_\lambda(x_i)-1
\bigr).
\]
Since $\BE_\rho L_{\lambda,\rho}(\Phi)=1$ we obtain that
%
\begin{equation}
\label{79} \BE_\rho D^n_{x_1,\ldots,x_n} L_{\lambda,\rho}(
\Phi) =\prod^n_{i=1}\bigl(h_\lambda(x_i)-1
\bigr),\qquad x_1,\ldots,x_n\in\BX, n\in\N.
\end{equation}

Denoting the right-hand side of (\ref{main1}) by $I$, we have
%
\begin{equation}
\label{81} I=\sum^\infty_{n=0}
\frac{1}{n!} \int\bigl(\BE_{\rho} L_{\lambda,\rho}(\Phi)
D^n_{x_1,\ldots,x_n}f(\Phi)\bigr) (\nu-\lambda)^n
\bigl(\mrmd(x_1,\ldots,x_n)\bigr).
\end{equation}
In the following, we assume that $f$ is bounded,
an assumption that will be removed in the final part
of the proof. Then $D^n_{x_1,\ldots,x_n}f(\Phi)$ is for
any fixed $(x_1,\ldots,x_n)$ bounded and hence square-integrable.
Hence, we can apply (\ref{28}) to the expectations in (\ref{81})
and use (\ref{79}) to obtain that
\begin{eqnarray*}
I&=&\sum^\infty_{n=0}\sum
^\infty_{k=0}\frac{1}{n!}\frac{1}{k!} \iint
\bigl(\BE_{\rho} D^{n+k}_{y_1,\ldots,y_k,x_1,\ldots,x_n}f(\Phi)\bigr) \prod
^k_{j=1}\bigl(h_\lambda(y_j)-1
\bigr)
\\
&&\hspace*{70.5pt}{}\times \rho^k\bigl(\mrmd(y_1,\ldots,y_k)\bigr) (\nu-
\lambda)^n\bigl(\mrmd(x_1,\ldots,x_n)\bigr)
\\
&=&\sum^\infty_{n=0}\sum
^\infty_{k=0}\frac{1}{n!}\frac{1}{k!} \int
\bigl(\BE_{\rho} D^{n+k}_{x_1,\ldots,x_{n+k}}f(\Phi)\bigr) \prod
^k_{j=1}\bigl(h_\lambda(x_j)-1
\bigr)
\\
&&\hspace*{62pt} {}\times\prod^{n+k}_{j=k+1}\bigl(h_\nu(x_j)-h_\lambda(x_j)
\bigr) \rho^{n+k}\bigl(\mrmd(x_1,\ldots,x_{n+k})\bigr),
\end{eqnarray*}
where the use of Fubini's theorem will be justified below.
Swapping the order of summation, we obtain
\begin{eqnarray*}
I&=&\sum^\infty_{k=0}\sum
^\infty_{n=k}\frac{1}{(n-k)!}\frac{1}{k!} \int
\bigl(\BE_{\rho} D^{n}_{x_1,\ldots,x_n}f(\Phi)\bigr)\prod
^k_{j=1}\bigl(h_\lambda(x_j)-1
\bigr)
\\
&&\hspace*{86pt} {}\times\prod^n_{j=k+1}\bigl(h_\nu(x_j)-h_\lambda(x_j)
\bigr)\rho^n\bigl(\mrmd(x_1,\ldots,x_n)\bigr)
\\
&=&\sum^\infty_{n=0}\frac{1}{n!}\int\sum
^n_{k=0}\pmatrix{n
\cr
k} \prod
^k_{j=1} \bigl(h_\lambda(x_j)-1
\bigr)\prod^n_{j=k+1}\bigl(h_\nu
(x_j)-h_\lambda(x_j)\bigr)
\\
&&\hspace*{37pt} {}\times\bigl(\BE_{\rho} D^{n}_{x_1,\ldots,x_n}f(\Phi)\bigr)
\rho^n\bigl(\mrmd(x_1,\ldots,x_n)\bigr)
\\
&=&\sum^\infty_{n=0}\frac{1}{n!} \int
\prod^n_{j=1}\bigl(h_\nu(x_j)-1
\bigr) \bigl(\BE_{\rho} D^{n}_{x_1,\ldots,x_n}f(\Phi)\bigr)
\rho^n\bigl(\mrmd(x_1,\ldots,x_n)\bigr),
\end{eqnarray*}
where we have used
\[
\prod^n_{j=1}\bigl(h_\nu(x_j)-h_\lambda(x_j)+h_\lambda(x_j)-1
\bigr)=\sum_{J\subset\{1,\ldots,n\}} \prod_{j_1\in J}
\bigl(h_\lambda(x_{j_1})-1\bigr)\prod
_{j_2\notin J}\bigl(h_\nu (x_{j_2})-h_\lambda(x_{j_2})
\bigr)
\]
and the permutation invariance of
$(\BE_{\rho} D^{n}_{x_1,\ldots,x_n}f(\Phi))\rho^n(\mrmd(x_1,\ldots,x_n))$
to obtain the last equality.
We are now using Proposition \ref{pbrown},
the identity (\ref{79}) with $\lambda$ replaced by $\nu$, and (\ref{28})
to obtain that
\begin{eqnarray*}
I&=&\sum^\infty_{n=0}\frac{1}{n!} \int
\bigl(\BE_{\rho}D^nL_{\nu,\rho}(\Phi)\bigr) \bigl(
\BE_{\rho}D^n f(\Phi)\bigr)\,\mrmd \rho^n
\\
&=&\BE_{\rho} L_{\nu,\rho}(\Phi)f(\Phi)=\BE_{\nu} f(\Phi),
\end{eqnarray*}
where $D^n f(\varphi)$ denotes for any
$\varphi\in\bN$ the mapping
$(x_1,\ldots,x_n)\mapsto D^n_{x_1,\ldots,x_n}f(\varphi)$.
This proves (\ref{main1}) for bounded $f$.

To justify the formal calculation above and to establish
(\ref{main2}), we need to show
that
\begin{eqnarray*}
c&:=&\sum^\infty_{n=0}\frac{1}{n!}\int
\sum^n_{k=0}\pmatrix{n
\cr
k} \prod
^k_{j=1}\bigl|h_\lambda(x_j)-1\bigr|\prod
^n_{j=k+1}\bigl|h_\nu(x_j)-h_\lambda
(x_j)\bigr|
\\
&&\hspace*{36pt} {}\times
\bigl|\BE_{\rho} D^{n}_{x_1,\ldots,x_n}f(\Phi)\bigr|\rho^n
\bigl(\mrmd(x_1,\ldots,x_n)\bigr)
\end{eqnarray*}
is finite. By permutation invariance,
\begin{eqnarray*}
c&=&\sum^\infty_{n=0}\frac{1}{n!}\int
\bigl(|h_\lambda-1|+|h_\nu-h_\lambda|\bigr)^{\otimes n} \bigl|
\BE_{\rho} D^{n}f(\Phi)\bigr|\,\mrmd \rho^n
\\
&\le&\sum^\infty_{n=0}\frac{1}{n!}\int
\bigl(2|h_\lambda-1|+|h_\nu-1|\bigr)^{\otimes n} \bigl|\BE_{\rho}
D^{n}f(\Phi)\bigr|\,\mrmd \rho^n.
\end{eqnarray*}
The Cauchy--Schwarz inequality yields,
\begin{eqnarray*}
c&\le&\sum^\infty_{n=0}\frac{\sqrt{a_n}}{n!}
\biggl(\int\bigl(\bigl(2|h_\lambda-1|+|h_\nu-1|\bigr)^{\otimes n}
\bigr)^2 \,\mrmd \rho^n \biggr)^{1/2}
\\
&=&\sum^\infty_{n=0}\frac{\sqrt{a_n}}{n!} \biggl(
\int\bigl(2|h_\lambda-1|+|h_\nu-1|\bigr)^2 \,\mrmd \rho
\biggr)^{n/2},
\end{eqnarray*}
where
\[
a_n:=\int\bigl(\BE_{\rho} D^{n}f(\Phi)
\bigr)^2 \,\mrmd \rho^n,\qquad n\in\N_0.
\]
Applying Cauchy--Schwarz again, yields
\[
c^2\le \Biggl(\sum^\infty_{n=0}
\frac{a_n}{n!} \Biggr) \sum^\infty_{n=0}
\frac{1}{n!} \biggl(\int\bigl(2|h_\lambda-1|+|h_\nu-1|\bigr)^2
\,\mrmd \rho \biggr)^{n}.
\]
The first series in the above product converges by (\ref{28})
(we have $\BE_{\rho} f(\Phi)^2<\infty$). The second
series converges, since the integral there is finite
by (\ref{22a}) and the Minkowski inequality.

We now extend the result to general $f$ satisfying
$\BE_\rho f(\Phi)^2<\infty$. We take a sequence
of bounded functions $f_l$, $l\in\N$, such that
$\BE_\rho(f(\Phi)-f_l(\Phi))^2\to0$ as $l\to\infty$.
We know already that
%
\begin{equation}
\label{1} \BE_{\nu} f_l(\Phi)= \sum
^\infty_{n=0}\frac{1}{n!}\int\bigl(
\BE_\lambda D^n f_l(\Phi)\bigr)\,\mrmd (\nu -
\lambda)^n
\end{equation}
holds for all $l\in\N$. By Cauchy--Schwarz,
%
\begin{eqnarray}
\label{1a} \BE_\nu\bigl|f(\Phi)-f_l(\Phi)\bigr|&=&
\BE_\rho L_{\nu,\rho}\bigl|f(\Phi )-f_l(\Phi)\bigr| \nonumber\\[-8pt]\\[-8pt]
&\le&\bigl(
\BE_\rho L_{\nu,\rho}^2\bigr)^{1/2} \bigl(
\BE_\rho\bigl(f(\Phi)-f_l(\Phi)\bigr)^2
\bigr)^{1/2} \to0\nonumber
\end{eqnarray}
as $l\to\infty$. Hence, the left-hand side of (\ref{1}) tends
to $\BE_\nu f(\Phi)$ as $l\to\infty$. To deal with
the right-hand side, we
consider sequences $\mathbf{g}=(g_n)_{n\ge0}$, where $g_0\in\R$
and $g_n$, $n\ge1$, is a measurable function on $\BX^n$.
Introduce the space $\mathbf{V}$ of all such sequences
satisfying
\[
\|\bg\|:=\sum^\infty_{n=0}\frac{1}{n!}
\int|g_n| |h_\nu-h_\lambda|^{\otimes n} \,\mrmd
\rho^n<\infty.
\]
Then $\mathbf{V}$ is a direct sum of Banach spaces and
hence a Banach space as well. For $l\in\N$ define
\[
g_{l,n}:=\BE_\lambda D^n f_l(\Phi),\qquad n
\ge0,\qquad \bg_l:=(g_{l,n})_{n\ge0}\in\mathbf{V}.
\]
Our next aim is to show that $(\bg_l)$ is a Cauchy-sequence.
We have for $l,m\in\N$ that
\begin{eqnarray*}
\|\bg_l-\bg_m\|&=& \sum^\infty_{n=0}
\frac{1}{n!} \int\bigl|\BE_\lambda D^n f_l(\Phi)-
\BE_\lambda D^n f_m(\Phi)\bigr| |h_\nu-h_\lambda|^{\otimes n}
\,\mrmd \rho^n
\\
&=&\sum^\infty_{n=0}\frac{1}{n!}\int\bigl|
\BE_\rho L_{\lambda,\rho
}(\Phi) D^n f_{l,m}(
\Phi)\bigr||h_\nu-h_\lambda|^{\otimes n} \,\mrmd \rho^n,
\end{eqnarray*}
where $f_{l,m}:=f_l-f_m$. From the calculation in the
first part of the proof, we obtain that
\[
\|\bg_l-\bg_m\|\le \sum^\infty_{n=0}
\frac{1}{n!}\int\bigl|\BE_\rho D^n f_{l,m}(\Phi)\bigr|
\bigl(2|h_\lambda-1|+|h_\nu-1|\bigr)^{\otimes n}\,\mrmd \rho^n.
\]
Applying the Cauchy--Schwarz inequality twice, as in the
second part of the proof yields
\[
\|\bg_l-\bg_m\|^2\le a\sum
^\infty_{n=0}\frac{1}{n!}\int\bigl(
\BE_\rho D^n f_{l,m}(\Phi)\bigr)^2 \,\mrmd
\rho^n,
\]
where
\[
a:=\sum^\infty_{n=0}\frac{1}{n!} \biggl(
\int\bigl(2|h_\lambda-1|+|h_\nu-1|\bigr)^2 \,\mrmd \rho
\biggr)^{n}.
\]
By (\ref{28}),
\[
\|\bg_l-\bg_m\|^2\le a \BE_\rho
f_{l,m}(\Phi)^2 =a \BE_\rho\bigl(f_l(
\Phi)-f_m(\Phi)\bigr)^2.
\]
By the choice of $f_l$ the sequence $(\bg_l)$ has
the Cauchy property. Because $\mathbf{V}$ is complete,
there is a $\bg=(g_n)\in\mathbf{V}$ such that
$\|\bg_l-\bg\|\to0$ as $l\to\infty$. Since,
\[
\Biggl|\sum^\infty_{n=0}\frac{1}{n!}\int
g_n\mathrm{d}(\nu-\lambda)^n -\sum^\infty_{n=0}
\frac{1}{n!}\int g_{l,n}\,\mrmd (\nu-\lambda)^n \Biggr| \le\sum
^\infty_{n=0}\frac{1}{n!}
\int|g_n-g_{l,n}| |h_\nu-h_\lambda|^{\otimes n}\,\mrmd
\rho^n,
\]
we obtain from (\ref{1}) and (\ref{1a}) that
\[
\BE_{\nu} f(\Phi)= \sum^\infty_{n=0}
\frac{1}{n!}\int g_n \mrmd(\nu-\lambda)^n.
\]
It remains to show that, for any fixed $n\ge0$,
%
\begin{equation}
\label{key} |h_\nu-h_\lambda|^{\otimes n}g_n=
|h_\nu-h_\lambda|^{\otimes n}\BE_\lambda
D^n f(\Phi),\qquad \rho^n\mbox{-a.e.}
\end{equation}
We claim that
%
\begin{equation}
\label{013} \lim_{l\to\infty}\int_{B^n}
\BE_\lambda\bigl| D^n f(\Phi) - D^n f_l(
\Phi)\bigr| \,\mrmd \rho^n=0
\end{equation}
for all $B\in\mathcal{X}$ with $\lambda(B)<\infty$ and
$\rho(B)<\infty$. As in the proof of \cite{LaPe11}, Lemma 2.3,
it suffices to demonstrate that
%
\begin{equation}
\label{014} \lim_{l\to\infty}\int_{B^m}
\BE_\lambda \Biggl|f \Biggl(\Phi+ \sum_{i=1}^m
\delta_{y_i} \Biggr)- f_l \Biggl(\Phi+ \sum
_{i=1}^m \delta_{y_i} \Biggr) \Biggr|
\rho^m\bigl(\mrmd(y_1,\ldots,y_m)\bigr)=0
\end{equation}
for all $m\in\{1,\ldots,n\}$.
By the (multivariate) Mecke equation
(see, e.g., \cite{Mecke} or \cite{LaPe11}, (2.10)) the integral in
(\ref{014})
equals
%
\begin{equation}
\label{456} \BE_\rho\int_{B^m} L_{\lambda,\rho}(
\Phi-\delta_{y_1}-\cdots-\delta_{y_m}) \bigl|f(
\Phi)-f_l(\Phi)\bigr|\Phi^{(m)}\bigl(\mrmd(y_1,\ldots,y_m)\bigr),
\end{equation}
where, for $\varphi\in\bN$,
$\varphi^{(m)}$ is the measure on $\BX^m$ defined by
%
\begin{eqnarray}
\label{phim} \varphi^{(m)}(C)&:=&\int\cdots\int\I_C(y_1,\ldots,y_m) \Biggl(\varphi- \sum_{j=1}^{m-1}
\delta_{y_j} \Biggr) (\mrmd y_m)
\nonumber\\
&&\hspace*{32pt}{}\times \Biggl(\varphi- \sum_{j=1}^{m-2}
\delta_{y_j} \Biggr) (\mrmd y_{m-1}) \times\cdots\\
&&\hspace*{32pt}{}\times (\varphi-
\delta_{y_1}) (\mrmd y_2)\varphi(\mrmd y_1),\qquad C\in
\mathcal{X}^{\otimes m}.\nonumber
\end{eqnarray}
By Lemma \ref{lkey} below and the Cauchy--Schwarz inequality,
(\ref{456}) tends to $0$ as $l\to\infty$.
Now (\ref{013}) implies that
$g_{l,n}=\BE_\lambda D^n f_l(\Phi)$ tends to
$\BE_\lambda D^n f(\Phi)$ $\rho^n$-a.e. on $B^n$ as $l\to\infty$
along a subsequence. Since
\[
\lim_{l\to\infty} \int|g_n-g_{l,n}||h_\nu-h_\lambda|^{\otimes n}\,\mrmd
\rho^n=0,
\]
there is a further subsequence $\N'\subset\N$ such that
$|h_\nu-h_\lambda|^{\otimes n}g_{l,n}$ tends to
$|h_\nu-h_\lambda|^{\otimes n}g_{n}$ $\rho^n$-a.e. on $B^n$
as $l\to\infty$ along $\N'$. It follows that
(\ref{key}) holds for $\rho^n$ restricted to $B^n$.
Since $\rho$ and $\lambda$ are $\sigma$-finite we obtain (\ref{key}).
This completes the proof of the theorem.
\end{pf*}

In the final part of the above proof, we have used
the following lemma. Recall the definition (\ref{phim}).
%
\begin{lemma}\label{lkey} Assume that (\ref{22a}) holds
and let $B\in\mathcal{X}$ satisfy $\lambda(B)<\infty$
and $\rho(B)<\infty$.
Then we have for all $m\ge1$ that
\[
\BE_\rho \biggl(\int_{B^m} L_{\lambda,\rho}(\Phi-
\delta_{x_1}-\cdots-\delta_{x_m}) \Phi^{(m)}
\bigl(\mrmd(x_1,\ldots,x_m)\bigr) \biggr)^2<
\infty.
\]
\end{lemma}
\begin{pf}
Writing the square of the inner integral as a double
integral and using a combinatorial argument,
we see that it suffices to prove that
\begin{eqnarray*}
&&
\BE_\rho\int_{B^{m-k}}\int_{B^m}
L_{\lambda,\rho}(\Phi-\delta_{x_1}-\cdots-\delta_{x_m})
L_{\lambda,\rho}(\Phi-\delta_{x_1}-\cdots-\delta_{x_k} -
\delta_{y_1}-\cdots-\delta_{y_{m-k}})
\\
&&\qquad\hspace*{34.3pt}{}\times(\Phi-\delta_{y_1}-\cdots-\delta_{y_{m-k}})^{(m)}
\bigl(\mrmd(x_1,\ldots,x_m)\bigr) \Phi^{(m-k)}
\bigl(\mrmd(y_1,\ldots,y_{m-k})\bigr)<\infty
\end{eqnarray*}
for all $k\in\{0,\ldots,m\}$ (with the obvious convention
for $k=m$).
Applying the Mecke equation twice,
we obtain that this expression equals
\begin{eqnarray*}
&&
\BE_\rho\int_{B^{m-k}}\int_{B^m}
L_{\lambda,\rho}(\Phi+\delta_{y_1}+\cdots+\delta_{y_{m-k}})
L_{\lambda,\rho}(\Phi+\delta_{x_{k+1}}+\cdots+\delta_{x_m})
\\
&&\qquad\hspace*{35.5pt}{}\times\rho^m\bigl(\mrmd(x_1,\ldots,x_m)\bigr)
\rho^{m-k}\bigl(\mrmd(y_1,\ldots,y_{m-k})\bigr).
\end{eqnarray*}
Since
\[
L_{\lambda,\rho}(\Phi+\delta_{y_1}+\cdots+\delta_{y_{m-k}})
=L_{\lambda,\rho}(\Phi)h_\lambda(y_1)\times\cdots\times
h_\lambda(y_{m-k}),
\]
we obtain that the above expectation equals
$
\rho(B)^{m-k}\lambda(B)^m\BE_\rho L_{\lambda,\rho}(\Phi)^2
$
which is finite by Proposition \ref{pbrown}.
\end{pf}
%
\begin{remark}\label{reps}
In the case $\rho=\lambda$ (this requires $\nu\ll\lambda$)
the proof of Theorem \ref{tvar3} becomes considerably
simpler. Another simplification is possible if
$\BE_\rho f(\Phi)^{2+\varepsilon}<\infty$ for
some $\varepsilon>0$.
Then $\BE_{\rho}(D^n_{x_1,\ldots,x_n} f(\Phi))^2<\infty$
for all $n\ge1$ and $\rho^n$-a.e. $(x_1,\ldots,x_n)$.
Indeed, by the proof of Lemma~2.3
in \cite{LaPe11} it is enough to show that
$
\BE_\rho f(\Phi)^2\Phi(B)^k<\infty
$
for all $k\in\N$ and any $B\in\mathcal{X}$ with $\rho(B)<\infty$.
Since $\Phi(B)$ has finite moments of any order, this
is a direct consequence of H\"older's inequality.
We can then apply (\ref{28}) to the expectations in (\ref{81})
and proceed exactly as in the proof of Theorem \ref{tvar3}.
This makes the final (and somewhat tricky) part of this
proof superfluous.
\end{remark}

We continue with providing special cases of Theorem \ref{tvar3}.
We let $\nu=\nu_1+\nu_2$ (resp., $\lambda=\lambda_1+\lambda_2$)
be the
Lebesgue decomposition of $\nu$ (resp., $\lambda$) with respect to
$\lambda$
(resp., $\nu$).
Hence $\nu_1\ll\lambda$ and $\nu_2\perp\lambda$,
where the latter means that $\nu_2$ and $\lambda$ are
\textit{singular}, that is concentrated on
disjoint measurable subsets of $\BX$.
%
\begin{theorem}\label{tvarleb2}
Let $f\dvtx \bN\rightarrow\R$ be measurable.
Assume that either
%
\begin{equation}
\label{44} \int \biggl(1-\frac{\mrmd \nu_1}{\mrmd \lambda} \biggr)^2 \,\mrmd \lambda+
\nu_2(\BX )<\infty
\end{equation}
and $\BE_{\lambda+\nu_2}f(\Phi)^2<\infty$, or that
%
\begin{equation}
\label{33} \int \biggl(1-\frac{\mrmd \lambda_1}{\mrmd \nu} \biggr)^2 \,\mrmd \nu+
\lambda_2(\BX )<\infty
\end{equation}
and $\BE_{\nu+\lambda_2}f(\Phi)^2<\infty$. Then (\ref{main1}) holds.
\end{theorem}
\begin{pf} We prove only the first assertion.
There are disjoint measurable subsets $B_1$ and
$B_2$ of $\BX$ such that
%
\begin{equation}
\label{A1} \lambda(\BX\setminus B_1)=\nu_2(\BX\setminus
B_2)=0.
\end{equation}
In particular, $\nu_1(\BX\setminus B_1)=0$.
Let $\rho:=\lambda+\nu_2$.
It is easy to check that
\[
h_\lambda=\I_{B_1},\qquad h_\nu=\I_{B_1}h_1+
\I_{B_2},
\]
where $h_1:=\mrmd \nu_1/\mrmd \lambda$. We have
\[
\int(h_\nu-1)^2 \,\mrmd \rho=\int(\I_{B_1}h_1-
\I_{\BX\setminus
B_2})^2\,\mrmd \rho =\int(\I_{B_1}h_1-
\I_{\BX\setminus B_2})^2\,\mrmd \lambda =\int(h_1-1)^2 \,\mrmd
\lambda
\]
and $\int(h_\lambda-1)^2 \,\mrmd \rho=\rho(\BX\setminus B_1)=\nu_2(\BX)$.
Therefore, (\ref{22a}) holds
and the result follows from Theorem \ref{tvar3}.
\end{pf}

The next corollary deals with a monotone perturbation
of $\lambda$.
%
\begin{corollary}\label{c1} Let $\mu$ be a $\sigma$-finite
measure on $\BX$ and assume that $h:=\mrmd \lambda/\mrmd(\lambda+\mu)$ satisfies
%
\begin{equation}
\label{square1} \int(1-h)^2\,\mrmd (\lambda+\mu)<\infty.
\end{equation}
Then we have for all measurable $f$ with
$\BE_{\lambda+\mu}f(\Phi)^2<\infty$ that
%
\begin{equation}
\label{m1} \BE_{\lambda+\mu} f(\Phi)=\BE_\lambda f(\Phi) +\sum
^\infty_{n=1}\frac{1}{n!}\int\bigl(
\BE_\lambda D^n_{x_1,\ldots,x_n} f(\Phi)\bigr)\mu^n
\bigl(\mrmd(x_1,\ldots,x_n)\bigr).
\end{equation}
\end{corollary}
\begin{pf} Apply the second part of Theorem \ref{tvarleb2} with $\nu
=\lambda+\mu$.
Then $\lambda_2=0$ and \mbox{$\mrmd \lambda/\mrmd \nu=h$}.
\end{pf}
%
\begin{remark}\label{r33}
In the situation of Corollary \ref{c1}, we
may assume that $h\le1$. Then $1-h$ is a
density of $\mu$ with respect to $\lambda+\mu$,
so that $\int(1-h)^2\,\mrmd (\lambda+\mu)=\int(1-h)\,\mrmd \mu$.
In particular, $\mu(\BX)<\infty$ implies (\ref{square1}),
cf. Theorem \ref{tvar}.
\end{remark}

The results of this section can be extended so as
to cover additional randomization.
%
\begin{remark}\label{rrand}
Let $(\BY,\mathcal{Y})$ be a measurable space and $\eta\dvtx
\Omega\rightarrow\BY$ be a measurable mapping such that
$\BP_\lambda((\eta,\Phi)\in\cdot)=\BV\otimes \Pi_\lambda$ for all
$\sigma$-finite measures $\lambda$, where $\BV$ is a probability
measure on $(\BY,\mathcal{Y})$, not depending on $\lambda$. The
definition of the difference operator can be extended to measurable
functions $f\dvtx \BY\times\bN\rightarrow\R$ in the following natural
way. If $n\in\N$ and $x_1,\ldots,x_n\in\BX$ then
$D^n_{x_1,\ldots,x_n}f\dvtx \BY\times\bN\rightarrow\R$ is defined by
$D^n_{x_1,\ldots,x_n}f(y,\varphi):=D^n_{x_1,\ldots,x_n}f_y(\varphi)$,
where $f_y:=f(y,\cdot)$, $y\in\BY$. Assume now that $\lambda,\nu,\rho$
satisfy the assumptions of Theorem \ref{tvar3} and that
\mbox{$\BE_{\rho}f(\eta,\Phi)^2<\infty$}. We claim that (\ref{main2}) and
(\ref{main1}) hold when replacing $\Phi$ by $(\eta,\Phi)$. This implies
that all results of this section (as well as those of Section
\ref{secder}) remain valid with the obvious changes.

To verify the above claim we define, for any $\varphi\in\bN$,
$
\tilde f(\varphi):=\int f(y,\varphi)\BV(\mrmd y)$
and conclude from Jensen's inequality that
$\BE_\rho\tilde f(\Phi)^2\le\BE_\rho f(\eta,\Phi)^2<\infty$.
Hence, Theorem \ref{tvar3} applies and we need to show
for all $n\in\N$ that
\[
\BE_\lambda D^n_{x_1,\ldots,x_n}\tilde{f}(\Phi)=
\BE_\lambda D^n_{x_1,\ldots,x_n}f(\eta,\Phi), \qquad
\rho^n\mbox{-a.e. $(x_1,\ldots,x_n)$}.
\]
In view of (\ref{Dsymmetric0}) and Fubini's theorem
it is sufficient to show for all $m\ge0$ that
\[
\BE_\lambda\int\bigl|f(y,\Phi+\delta_{x_1}+\cdots+
\delta_{x_m})\bigr|\BV(\mrmd y)= \BE_\lambda\bigl|f(\eta,\Phi+
\delta_{x_1}+\cdots+\delta_{x_m})\bigr|<\infty
\]
for $\rho^m$-a.e. $(x_1,\ldots,x_m)$ (with the obvious convention for
$m=0$). To this end, we take $B_1,\ldots,B_m\in\mathcal{X}$
with finite measure with respect to both $\lambda$ and $\rho$,
let $B:=B_1\times\cdots\times B_m$, and obtain from the Mecke equation that
\begin{eqnarray*}
&&\int_B\BE_\lambda\bigl|f(\eta,\Phi+
\delta_{x_1}+\cdots+\delta_{x_m})\bigr|\rho^m
\bigl(\mrmd(x_1,\ldots,x_m)\bigr)
\\
&&\quad=\BE_\rho\int_B L_{\lambda,\rho}(\Phi) \bigl|f(
\eta,\Phi+\delta_{x_1}+\cdots+\delta_{x_m})\bigr|\rho^m
\bigl(\mrmd(x_1,\ldots ,x_m)\bigr)
\\
&&\quad=\BE_\rho\bigl|f(\eta,\Phi)\bigr|\int_B
L_{\lambda,\rho}(\Phi-\delta_{x_1}-\cdots-\delta_{x_m})
\Phi^{(m)}\bigl(\mrmd(x_1,\ldots,x_m)\bigr),
\end{eqnarray*}
which is finite by Cauchy--Schwarz,
Lemma \ref{lkey} and our assumption $\BE_{\rho}f(\eta,\Phi
)^2<\infty$.
\end{remark}

\section{Necessary conditions for the variational formulas}\label{secnec}

Again we consider two $\sigma$-finite measures $\lambda,\nu$ on $\BX$.
The squared \textit{Hellinger distance} between these two
measures is defined as
%
\begin{equation}
\label{Hell} H(\lambda,\nu):=\frac{1}{2}\int (\sqrt{h_\lambda}-
\sqrt {h_\nu} )^2 \,\mrmd \rho,
\end{equation}
where (as before) $\rho$ is a $\sigma$-finite measure dominating
$\lambda$ and $\nu$
and $h_\lambda$, respectively, $h_\nu$ are the corresponding densities.
%
\begin{theorem}\label{tnec} Assume that
(\ref{main1}) holds for all bounded measurable
$f\dvtx \bN\rightarrow\R$. Then $\Pi_\lambda$ and $\Pi_\nu$ are not singular
and
%
\begin{equation}
\label{370} H(\lambda,\nu)<\infty.
\end{equation}
\end{theorem}
\begin{pf} Assume on the contrary that $\Pi_\lambda$ and $\Pi_\nu
$ are singular.
Then we find disjoint sets $F,G\in\mathcal{N}$ such that
$\Pi_\lambda(F)=\Pi_{\nu}(G)=1$.
We now proceed as in the proof of Theorem 9.1.13 in \cite{MKM82}.
Let $C_n\in\mathcal{X}$, $n\in\N$, be such that $\lambda(C_n)+\nu
(C_n)<\infty$,
and $C_n\uparrow\BX$ as $n\to\infty$. Recall that the
restriction of $\varphi\in\bN$ to $B\in\mathcal{X}$ is
denoted by $\varphi_B$.
We have for any $n\in\N$ that
\begin{eqnarray*}
\exp\bigl[-\lambda(C_n)\bigr]&=&\BP_\lambda\bigl(
\Phi(C_n)=0,\Phi\in F\bigr) =\BP_\lambda\bigl(
\Phi(C_n)=0,\Phi_{\BX\setminus C_n}\in F\bigr)
\\
&=&\BP_\lambda\bigl(\Phi(C_n)=0\bigr)\BP_\lambda(
\Phi_{\BX\setminus C_n}\in F) =\exp\bigl[-\lambda(C_n)\bigr]
\BP_\lambda(\Phi_{\BX\setminus C_n}\in F).
\end{eqnarray*}
A similar calculation applies to $\BP_\nu$.
It follows that the sets
%
\begin{equation}
F_n:=\bigcap^\infty_{m=n}\{
\varphi\in\bN\dvtx \varphi_{\BX\setminus
C_n}\in F\},\qquad G_n:=\bigcap
^\infty_{m=n}\{\varphi\in\bN\dvtx \varphi_{\BX\setminus
C_n}\in
G\},\qquad n\in\N,
\end{equation}
have the properties
\[
\BP_\lambda(\Phi\in F_n)=\BP_{\nu}(\Phi\in
G_n)=1,\qquad n\in\N.
\]
This implies
%
\begin{equation}
\label{contra} \BP_\lambda\bigl(\Phi\in F'\bigr)=
\BP_{\nu}\bigl(\Phi\in G'\bigr)=1,
\end{equation}
where $F':=\bigcup_{n\in\N}F_n$ and $G':=\bigcup_{n\in\N}G_n$. Since
$F\cap G=\varnothing$ we have $F_n\cap G_n=\varnothing$ for any
$n\in\N$. Since $F_n$ and $G_n$ are increasing, we obtain that $F'\cap
G'=\varnothing$. Since $C_n\uparrow\BX$ we have for any
$(\varphi,x)\in\bN\times \BX$ that $\varphi\in F'$ if and only if
$\varphi+\delta_x\in F'$. Therefore,\vspace*{1pt} for $f:=\I_{F'}$,
$D^n_{x_1,\ldots,x_n}f\equiv0$ for all $n\in\N$ and all
$x_1,\ldots,x_n\in\BX$. Using this fact as well as (\ref{contra})
(together with $F'\cap G'=\varnothing$), we see that (\ref{main1})
fails.

A classical result by Liese \cite{Liese76} (see also \cite
{LieseVajda87}, Theorem
(3.30))
says that
%
\begin{equation}
\label{Hellinger}
H(\Pi_\lambda,\Pi_\nu)=1-\RMe^{-H(\lambda,\nu)}
\end{equation}
so that singularity of $\Pi_\lambda$ and $\Pi_\nu$ is
equivalent to $H(\lambda,\nu)=\infty$ (see
\cite{Liese76}, (3.2)).
\end{pf}

Recall the Lebesgue decompositions
$\nu=\nu_1+\nu_2$ of $\nu$ with respect to $\lambda$
and $\lambda=\lambda_1+\lambda_2$ of $\lambda$ with respect to $\nu$.
%
\begin{corollary}\label{tnec2} Assume that
(\ref{main1}) holds for all bounded measurable
$f\dvtx \bN\rightarrow\R$. Then
%
\begin{equation}
\label{37} \int \biggl(1-\sqrt{\frac{\mrmd \nu_1}{\mrmd \lambda}} \biggr)^2 \,\mrmd \lambda+
\nu_2(\BX) +\int \biggl(1-\sqrt{\frac{\mrmd \lambda_1}{\mrmd \nu}}
\biggr)^2 \,\mrmd \nu +\lambda_2(\BX)<\infty.
\end{equation}
Moreover, we have that $\Pi_{\nu_1}\ll\Pi_\lambda$
and $\Pi_{\lambda_1}\ll\Pi_\nu$ and in particular
$\Pi_{\nu}\ll\Pi_\lambda$ (resp., $\Pi_{\lambda}\ll\Pi_\nu$) provided
that $\lambda\ll\nu$ (resp., $\nu\ll\lambda$).
If, in addition, the density $\mrmd \nu_1/\mrmd \lambda$
(resp., $\mrmd \lambda_1/\mrmd \nu$) may be chosen bounded,
then (\ref{44}) (resp., (\ref{33})) holds.
\end{corollary}
\begin{pf} Since the definition (\ref{Hell}) of $H(\lambda,\nu)$
is independent of the dominating measure $\rho$, we have (see
also the proof of Theorem \ref{tvarleb2})
%
\begin{equation}
\label{H} H(\lambda,\nu)= \int \biggl(1-\sqrt{\frac{\mrmd \nu_1}{\mrmd \lambda}}
\biggr)^2 \,\mrmd \lambda+\nu_2(\BX) =\int \biggl(1-\sqrt{
\frac{\mrmd \lambda_1}{\mrmd \nu}} \biggr)^2 \,\mrmd \nu +\lambda_2(\BX).
\end{equation}
Hence, (\ref{37}) follows from (\ref{370}) while the asserted
absolute continuity relations follow from (\ref{37})
and \cite{Liese76}, Satz (3.3) (see \cite{MKM82}, Theorem 1.5.12).
If $\mrmd \nu_1/\mrmd \lambda$ may be chosen bounded,
then (\ref{33}) follows from (\ref{37}) and
the identity $(1-x)=(1-\sqrt{x})(1+\sqrt{x})$,
$x\ge0$.
\end{pf}

For monotone perturbations,
Corollary \ref{tnec2} yields the following characterization
of the variational formula.
%
\begin{corollary}\label{c4} Let $\mu$ be a $\sigma$-finite
measure on $\BX$.
\begin{longlist}[(ii)]
\item[(i)]
The variational formula (\ref{m1})
holds for all bounded and measurable $f\dvtx \bN\rightarrow\R$
if and only if $h:=\mrmd \lambda/\mrmd(\lambda+\mu)$ satisfies (\ref{square1}).
\item[(ii)]
Assume that $\mu\le\lambda$ and let $\nu:=\lambda-\mu$. Then
(\ref{main1}) holds for all bounded and measurable $f\dvtx \bN\rightarrow
\R$
if and only if $h_\mu:=\mrmd \mu/\mrmd \lambda$ satisfies $\int h_\mu^2\,\mrmd \lambda<\infty$.
\end{longlist}
\end{corollary}
%
\begin{remark}\label{ralt}
In general, inequality (\ref{37}) is weaker than both (\ref{44})
and (\ref{33}). We do not know whether
(\ref{37}) is sufficient for (\ref{main1})
to hold for all bounded measurable $f$.
\end{remark}
%
\begin{example}\label{excounter}
Assume that $\lambda$ is Lebesgue measure on $\BX:=\R^d$ for some
$d\ge1$. Let $\mu:=c\lambda$ for some $c>0$. Then
$\mrmd\lambda/\mrmd(\lambda+\mu)=(1+c)^{-1}$, so that (\ref {square1}) fails.
Let $B_n$ be a ball with centre at the origin and radius $n\in\N$ and
let $f$ be the measurable function on $\bN$ defined by
\[
f(\varphi):=\I\Bigl\{\lim_{n\to\infty}\lambda(B_n)^{-1}
\varphi (B_n)=1\Bigr\}.
\]
Then $\BE_\lambda f(\Phi)=1$ while
$\BE_{\lambda+\mu} f(\Phi)=0$. On the other hand we have
$D^n_{x_1,\ldots,x_n}f\equiv0$ for all $n\ge1$ and
all $x_1,\ldots,x_n\in\R^d$. Hence (\ref{m1}) fails.
\end{example}
%
\begin{remark}\label{remcomp}
Theorem \ref{tnec} and (\ref{H}) show that
(\ref{main1}) can only hold for all bounded functions
$f$ if the
 non-absolutely continuous part of the perturbation of $\lambda$
has finite mass while the
absolutely continuous part of the perturbation
leads to a distribution $\Pi_\nu$ that is absolutely continuous
with respect to the original distribution $\Pi_\lambda$.
Example \ref{excounter} shows what can go wrong with (\ref{main1})
if this second condition fails.
If one condition is violated, then this does not mean that
(\ref{main1}) does not hold for \textit{some} bounded measurable $f$.
In fact, Theorem \ref{tvar3} shows that the formula holds whenever $f$ depends
on the restriction of $\Phi$ to a set $B\in\mathcal{X}$ with
$\lambda(B)<\infty$ and $\nu(B)<\infty$.
\end{remark}

\section{Derivatives and Russo-type formulas}\label{secder}

In this section, we consider $\sigma$-finite measures
$\lambda,\rho$ on $\BX$
and assume that $\lambda$ is absolutely continuous with
respect to $\rho$ with density $h_\lambda$.
We also consider a measurable function $h\dvtx \BX\rightarrow\R$
and assume that
%
\begin{equation}
\label{56} \int(1-h_\lambda)^2\,\mrmd \rho+\int h^2\,\mrmd
\rho<\infty.
\end{equation}

\begin{theorem}\label{tseries} Assume that (\ref{56}) holds.
Let $\theta_0\in\R$ and assume that $I\subset\R$ is an interval
with non-empty interior such that $\theta_0\in I$ and
$h_\theta:=h_\lambda+(\theta-\theta_0)h\ge0$ $\rho$-a.e. for
$\theta\in I$.
For $\theta\in I$ let $\lambda_\theta$ denote the measure
with density $h_\theta$ with respect to $\rho$.
Let $f\dvtx \bN\rightarrow\R$ be a measurable function such that
$\BE_{\rho}f(\Phi)^2<\infty$.
Then,
%
\begin{equation}
\label{ps} \BE_{\lambda_\theta} f(\Phi)=\BE_\lambda f(\Phi) +\sum
^\infty_{n=1}\frac{(\theta-\theta_0)^n}{n!}\int\bigl(
\BE_\lambda D^nf(\Phi)\bigr)h^{\otimes n}\,\mrmd
\rho^n,\qquad \theta\in I,
\end{equation}
where $\BE_\lambda D^n(\Phi)$ denotes the function
$(x_1,\ldots,x_n)\mapsto\BE_\lambda D^n_{x_1,\ldots,x_n}f(\Phi)$
and the series converges absolutely.
Moreover,
%
\begin{equation}
\label{deriv} \frac{\mrmd}{\mrmd \theta}\BE_{\lambda_\theta} f(\Phi) =\int\bigl(
\BE_{\lambda_\theta} D_xf(\Phi)\bigr)h(x)\rho(\mrmd x), \qquad\theta \in I.
\end{equation}
\end{theorem}
\begin{pf} Let $\theta\in I$. By our assumptions
$1-h_\theta=(1-h_\lambda)-(\theta-\theta_0)h$
is square-integrable with respect to $\rho$. Hence
we can apply Theorem \ref{tvar3} with $\nu=\lambda_\theta$ to obtain
(\ref{ps}). In particular we get (\ref{deriv}) for $\theta=\theta_0$.

To derive (\ref{deriv}) for general $\theta\in I$ we apply the above
with $(\lambda_\theta,h_\theta)$ instead of $(\lambda,h_\lambda)$
and with $\theta$ instead of $\theta_0$.
Since
\[
h_{\tilde\theta}=h_\lambda+(\theta-\theta_0)h +(\tilde
\theta-\theta)h=h_\theta+(\tilde\theta-\theta)h, \qquad\tilde\theta\in I,
\]
we obtain the desired result from (\ref{deriv}) using the
same function $h$ as before.
\end{pf}
%
\begin{corollary}\label{c43} Let $\nu$ be another $\sigma$-finite
measure with
density $h_\nu$ with respect to $\rho$. Assume that (\ref{22a}) holds.
Then
%
\begin{equation}
\label{ps2} \BE_{\lambda+\theta(\nu-\lambda)} f(\Phi)=\BE_\lambda f(\Phi) +\sum
^\infty_{n=1}\frac{\theta^n}{n!}\int\bigl(
\BE_\lambda D^nf(\Phi )\bigr)\,\mrmd (\nu-\lambda)^n,\qquad
\theta\in[0,1],
\end{equation}
provided that $\BE_{\rho}f(\Phi)^2<\infty$.
\end{corollary}
\begin{pf} We take in Theorem \ref{tseries} $h:=h_\nu-h_\lambda$,
$I:=[0,1]$ and $\theta_0:=0$. The result follows
upon noting that square-integrability of $h$ is implied by
the Minkowski inequality.
\end{pf}
%
\begin{remark}\label{rgat}
Fix a measurable function $f\dvtx \bN \rightarrow\R$ such that
$\BE_{\rho}f(\Phi)^2<\infty$. Let $h_\lambda$ satisfy
$\int(1-h_\lambda)^2\,\mrmd \rho<\infty$ and let $H_\lambda$ be the set of all
measurable functions $h\dvtx \BX\rightarrow\R$ such that $\int
h^2\,\mrmd \rho<\infty$ and $h_\lambda+\theta h\ge0$ $\rho$-a.e. for all
$\theta$ in some (possibly one-sided) neighborhood $I_h$ of $0$. For
$h\in H_\lambda$ and $\theta\in I_h$ we let $\lambda_\theta$ denote the
measure with density $h_\theta:=h_\lambda+ \theta h$ with respect to
$\rho$. Then Theorem \ref{tseries} states that
%
\begin{equation}
\label{865} \lim_{\theta\to0} \theta^{-1}\bigl(\BE_{\lambda_\theta}
f(\Phi)-\BE_{\lambda} f(\Phi)\bigr) =G_{\lambda,f}(h),\qquad h\in
H_\lambda,
\end{equation}
where
%
\begin{equation}
\label{Glambda} G_{\lambda,f}(h):=\int\bigl(\BE_{\lambda}
D_xf(\Phi)\bigr)h(x)\rho(\mrmd x).
\end{equation}
Hence $G_{\lambda,f}(h)$ is the \textit{G\^ateaux derivative} of
the mapping $\nu\mapsto\BE_{\nu}f(\Phi)$
at $\lambda$ in the \textit{direction}~$h$.
\end{remark}

If the perturbation is absolutely continuous with respect to
the original measure $\lambda$, then we can strengthen (\ref{865})
to \textit{Fr\'echet} differentiability as follows.
Let $H^*_\lambda$ be the set
of all measurable functions $h\dvtx \BX\rightarrow\R$ such that
$\int h^2\,\mrmd \lambda<\infty$ and $1 +h\ge0$ $\lambda$-a.e.
%
\begin{proposition}\label{rfrech} Let $f\dvtx \bN\rightarrow\R$ be measurable
and such that $\BE_{\lambda}f(\Phi)^2<\infty$.
For $h\in H^*_\lambda$ let $\lambda_h$ denote the measure with density
$1+h$ with respect to $\lambda$. Then
%
\begin{equation}
\label{866} \BE_{\lambda_h} f(\Phi)=\BE_{\lambda} f(
\Phi)+G_{\lambda,f}(h) +\RMo\bigl(\|h\|\bigr),\qquad h\in H^*_\lambda,
\end{equation}
where $G_{\lambda,f}(h)$ is defined by (\ref{Glambda}),
$\|h\|:=\sqrt{\int h^2 \,\mrmd \lambda}$ and
$\lim_{t\to0} t^{-1}\RMo(t)=0$.
\end{proposition}
\begin{pf}
We apply Theorem \ref{tvar3} with $\rho=\lambda$
(so that $h_\lambda\equiv1$) and $\nu=\lambda_h$ to obtain that
\[
\BE_{\lambda_h} f(\Phi)=\BE_{\lambda} f(\Phi)+G_{\lambda,f}(h)+c_h,
\]
where
\[
c_h:=\sum^\infty_{n=2}
\frac{1}{n!} \int\BE_\lambda D^nf(\Phi)h^{\otimes n}\,\mrmd
\lambda^n.
\]
Applying the triangle inequality and then the
Cauchy--Schwarz inequality to each summand
gives
\[
|c_h|\le\sum^\infty_{n=2}
\frac{1}{n!} \biggl(\int\bigl(\BE_\lambda D^nf(\Phi)
\bigr)^2\,\mrmd \lambda^n \biggr)^{1/2} \biggl(\int
h^2\,\mrmd \lambda \biggr)^{n/2}.
\]
Applying the Cauchy--Schwarz inequality again yields
\[
|c_h|\le \Biggl(\sum^\infty_{n=2}
\frac{1}{n!} \int\bigl(\BE_\lambda D^nf(\Phi)
\bigr)^2\,\mrmd \lambda^n \Biggr)^{1/2} \Biggl(\sum
^\infty_{n=2}\frac{1}{n!} \biggl(\int
h^2\,\mrmd \lambda \biggr)^{n} \Biggr)^{1/2}.
\]
The first factor is finite by (\ref{28})
and the second equals $\tilde{\mathrm{o}}(\|h\|)$, where
$\tilde{\mathrm{o}}(t):=\sqrt{\RMe^{t^2}-1-t^2}$.
\end{pf}

Next, we generalize (\ref{deriv}) to possibly non-linear
perturbations of $\lambda$.
%
\begin{theorem}\label{tderiv} Assume that (\ref{56}) holds.
Let $\theta_0\in\R$ and assume that $I\subset\R$ is an interval
with non-empty interior such that $\theta_0\in I$. For any $\theta\in
I$ let
$R_\theta\dvtx \BX\rightarrow\R$ be a measurable function
such that the following assumptions are satisfied:
\begin{longlist}[(iii)]
\item[(i)] For all $\theta\in I$, $h_\lambda+(\theta-\theta_0)(h+R_\theta)\ge0$
$\rho$-a.e.
\item[(ii)] $\lim_{\theta\to\theta_0}R_\theta= 0$ $\rho$-a.e.
\item[(iii)] There is a measurable function
$R\dvtx \BX\rightarrow[0,\infty)$ such that $|R_\theta|\le R$
$\rho$-a.e. for all $\theta\in I$ and $\int R^2\,\mrmd \rho<\infty$.
\end{longlist}
For $\theta\in I$, let $\lambda_\theta$ denote the measure
with density $h_\lambda+(\theta-\theta_0)(h+R_\theta)$ with respect
to $\rho$.
Let $f\dvtx \bN\rightarrow\R$ be a measurable function such that
$\BE_{\rho}f(\Phi)^2<\infty$. Then
%
\begin{equation}
\label{deriv5} \frac{\mrmd}{\mrmd \theta}\BE_{\lambda_\theta} f(\Phi) \bigg|_{\theta
=\theta_0} =
\int\bigl(\BE_{\lambda} D_xf(\Phi)\bigr)h(x)\rho(\mrmd x).
\end{equation}
\end{theorem}
\begin{pf} In view of $\int h^2\,\mrmd \rho<\infty$ and assumption (iii),
it is possible to apply Theorem \ref{tvar3} to the measure $\nu
=\lambda_\theta$.
This gives for $\theta\in I\setminus\{\theta_0\}$
%
\begin{eqnarray}
\label{ps5}
&&(\theta-\theta_0)^{-1}\bigl(
\BE_{\lambda_\theta} f(\Phi)-\BE_\lambda f(\Phi)\bigr)\nonumber\\
&&\quad=\int\bigl(
\BE_\lambda Df(\Phi)\bigr) (h+R_\theta) \,\mrmd \rho
\\
&&\qquad{}+\sum^\infty_{n=2}\frac{(\theta-\theta_0)^{n-1}}{n!} \int
\bigl(\BE_\lambda D^nf(\Phi)\bigr) (h+R_\theta)^{\otimes n}\,\mrmd
\rho^n.\nonumber
\end{eqnarray}
Applying Theorem \ref{tvar3} to the measure $\nu$ with density
$h_\lambda+|h|+R$ with respect to $\rho$ yields
\[
\sum^\infty_{n=0}\frac{1}{n!}\int\bigl|
\BE_\lambda D^nf(\Phi)\bigr| \bigl(|h|+R\bigr)^{\otimes n}\,\mrmd
\rho^n<\infty.
\]
Hence the result follows from assumption (ii) and bounded
convergence.
\end{pf}

The case where the perturbed measure $\lambda_\theta$
is absolutely continuous with respect to $\lambda$
is of special interest. Then the assumptions (ii) and (iii) in
Theorem \ref{tderiv} can be simplified.
%
\begin{theorem}\label{t57}
Assume that $\int h^2\,\mrmd \lambda<\infty$.
Let $\theta_0\in\R$ and assume that $I\subset\R$ is an interval
with non-empty interior such that $\theta_0\in I$. For any $\theta\in
I$ let
$R_\theta\dvtx \BX\rightarrow\R$ be a measurable function
such that the following assumptions are satisfied:
\begin{longlist}[(ii)]
\item[(i)] For all $\theta\in I$, $1+(\theta-\theta_0)(h+R_\theta)\ge0$ $\lambda$-a.e.
\item[(ii)] $\lim_{\theta\to\theta_0}\int R^2_\theta
\,\mrmd \lambda= 0$.
\end{longlist}
For $\theta\in I$, let $\lambda_\theta$ denote the measure
with density $1+(\theta-\theta_0)(h+R_\theta)$ with respect to
$\lambda$.
Let $f\dvtx \bN\rightarrow\R$ be a measurable function such that
$\BE_{\lambda}f(\Phi)^2<\infty$. Then
%
\begin{equation}
\label{deriv0} \frac{\mrmd}{\mrmd \theta}\BE_{\lambda_\theta} f(\Phi) \bigg|_{\theta
=\theta_0} =
\int\bigl(\BE_{\lambda} D_xf(\Phi)\bigr)h(x)\lambda(\mrmd x).
\end{equation}
\end{theorem}
\begin{pf} This time we apply Theorem \ref{tvar3} with $\rho
=\lambda$
(so that $h_\lambda\equiv1$) and $\nu=\lambda_\theta$. To treat
the right-hand side of (\ref{ps5}), we first note that
\[
\int\bigl|\BE_\lambda Df(\Phi)\bigr||R_\theta| \,\mrmd \lambda \le \biggl(\int
\bigl(\BE_\lambda Df(\Phi)\bigr)^2\,\mrmd \lambda
\biggr)^{1/2} \biggl(\int R_\theta^2 \,\mrmd \lambda
\biggr)^{1/2}.
\]
By assumption (ii), this tends to zero as $\theta\to\theta_0$. It
remains to
show that
\[
c_\theta:=\sum^\infty_{n=2}
\frac{1}{n!} \int\bigl|\BE_\lambda D^nf(\Phi)\bigr||h+R_\theta|^{\otimes n}\,\mrmd
\lambda^n
\]
is bounded in $\theta$. As in the proof of Proposition \ref{rfrech},
it follows
that
\[
c^2_\theta\le \Biggl(\sum^\infty_{n=2}
\frac{1}{n!}\int\bigl(\BE_\lambda D^nf(\Phi)
\bigr)^2\,\mrmd \lambda^n \Biggr) \Biggl(\sum
^\infty_{n=2}\frac{1}{n!} \biggl(
\int(h+R_\theta )^2\,\mrmd \lambda \biggr)^{n} \Biggr).
\]
Here the first factor is finite by Theorem \ref{tvar3}
while the second remains bounded by assumption~(ii).
\end{pf}
%
\begin{corollary}\label{t55}
Let the assumptions of Theorem \ref{t57} be satisfied.
Then
%
\begin{equation}
\label{deriv52} \frac{\mrmd}{\mrmd \theta}\BE_{\lambda_\theta} f(\Phi) \bigg|_{\theta
=\theta_0} =
\BE_{\lambda} \int\bigl(f(\Phi)-f(\Phi-\delta_x)\bigr)h(x)
\Phi(\mrmd x).
\end{equation}
\end{corollary}
\begin{pf} The result follows from (\ref{deriv0}) and
the Mecke equation from \cite{Mecke}.
\end{pf}
%
\begin{remark}\label{rlit}
The results of this section generalize the Poisson cases of the
derivative formulas in \cite{BaKlZu95} and \cite{Decreuse98}, where one
can also find some earlier predecessors. We note that \cite{BaKlZu95}
and \cite{Decreuse98} study more general point processes.
\end{remark}

Finally in this section, we deal with the case, where $\lambda_\theta
$ is
a multiple of a finite measure.

\begin{corollary}\label{finite} Assume that $\lambda$ is a finite measure
and let $f\dvtx \bN\rightarrow\R$ be a measurable function such that
$\BE_{\theta_0\lambda}f(\Phi)^2<\infty$ for some $\theta_0>0$.
Then $\theta\mapsto\BE_{\theta\lambda}f(\Phi)$
is analytic on $[0,\infty)$. Moreover,
%
\begin{eqnarray}
\label{deriv4} \frac{\mrmd}{\mrmd \theta}\BE_{\theta\lambda} f(\Phi) &=&
\int
\BE_{\theta\lambda} D_xf(\Phi)\lambda(\mrmd x),\qquad \theta\ge 0,
\\
\label{deriv09} \frac{\mrmd}{\mrmd \theta}\BE_{\theta\lambda} f(\Phi)&=&
\theta^{-1}\BE_{\theta\lambda}\int\bigl(f(\Phi)-f(\Phi-
\delta_x)\bigr) \Phi(\mrmd x),\qquad \theta> 0.
\end{eqnarray}
\end{corollary}
\begin{pf} Apply Theorem \ref{tseries}
with $\rho:=\theta_0\lambda$, $\lambda:=0$, $h:=\theta_0^{-1}$,
$\theta_0:=0$ and
$I:=[0,\infty)$.
This yields the first two assertions. As before, formula
(\ref{deriv09}) is a consequence of (\ref{deriv4}) and the Mecke
formula.
\end{pf}
%
\begin{remark}\label{r370} Consider in Corollary \ref{finite}
a general $\sigma$-finite measure $\lambda$ but assume
that the function $f$ does only depend on the restriction
of $\Phi$ to some set $B\in\mathcal{X}$ with $\lambda(B)<\infty$.
Applying the corollary to $\lambda(B\cap\cdot)$ gives (\ref{deriv4}).
This is extended in \cite{BordTo08} to
functions that depend measurably on the $\sigma$-field associated
with a \textit{stopping set} satisfying suitable integrability
assumptions.
\end{remark}
%
\begin{remark}\label{rrusso}
Let $f:=\I_A$, where $A\in\mathcal{N}$ is \textit{increasing}, that is,
whenever $\varphi\in A$ then $\varphi+\delta_x\in A$ for all $x\in\BX$.
Then
\[
\int\bigl(f(\Phi)-f(\Phi-\delta_x)\bigr)\Phi(\mrmd x) =\int\I\{\Phi\in
A,\Phi-\delta_x\notin A\}\Phi(\mrmd x)
\]
is the number of points of $\Phi$ that are \textit{pivotal} for $A$.
Hence (\ref{deriv09}) expresses the derivative of
$\BP_{\theta\lambda}(\Phi\in A)$ in terms
of the expected number of pivotal elements.
This Poisson counterpart
of the \textit{Margulis--Russo formula} for Bernoulli fields
was first proved in \cite{Zue92b}.
In the more general setting of Corollary \ref{t55},
the pivotal elements have to be counted in a weighted way.
\end{remark}

\section{Perturbation analysis of L\'evy processes}\label{seclevy}

In this section, we apply our results to $\R^d$-valued \textit{L\'evy
processes}, that is, to processes $X=(X_t)_{t\ge0}$ with homogeneous
and independent increments and $X_0=0$. We assume that $X$ is
continuous in probability. By Proposition II.3.36 in \cite{JaSh87} and
Theorem 15.4 in \cite{Kallenberg}, we can then assume that a.s.
%
\begin{equation}
\label{LCh} X_t=bt+W_t+\int_{|x|\le1}\int
^t_0x\bigl(\Phi(\mrmd s,\mrmd x)-\mrmd s\nu(\mrmd x)\bigr) +\int
_{|x|> 1}\int^t_0x \Phi(\mrmd s,\mrmd x),\qquad
t\ge0,\quad
\end{equation}
where $b\in\R^d$, $W=(W)_{t\ge0}$ is a $d$-dimensional Wiener process
with covariance matrix $\Sigma$ and $\Phi$ is an independent Poisson
process on $[0,\infty)\times\R^d$ with intensity measure
$\lambda_1\otimes\nu$. Here $\lambda_1$ is Lebesgue measure on
$[0,\infty)$ and $\nu$ is a \textit{L\'evy measure} on $\R^d$, that is,
a measure on $\R^d$ having $\nu(\{0\})=0$, and $ \int(|x|^2\wedge1)
\nu(\mrmd x)<\infty. $ The integrals in (\ref{LCh}) have to be interpreted
as limits in probability. Let $\bD$ denote the space of all
$\R^d$-valued right-continuous functions on $\R_+$ with left-hand
limits on $(0,\infty)$. By \cite{Kallenberg}, Theorem 15.1, we can and
will interpret $X$ as a random element in $\bD$ equipped with the
Kolmogorov product $\sigma$-field. The \textit{characteristic triplet}
$(\Sigma,b,\nu)$ determines the distribution of $X$. In this section,
we fix $\Sigma$ and let $\BP_{b,\nu}$ denote a probability measure on
$(\Omega,\mathcal{F})$ such that $\BP_{b,\nu}(X\in\cdot)$ is the
distribution of a L\'evy process with characteristic triplet
$(\Sigma,b,\nu)$. The expectation with respect to this measure is
denoted by $\BE_{b,\nu}$. As before, we let $\BP_{\lambda_1\otimes\nu}$
denote a probability measure such that
$\BP_{\lambda_1\otimes\nu}(\Phi\in\cdot)=\Pi_{\lambda _1\otimes\nu}$.
Similarly as in Remark~\ref{rrand}, we assume that under
$\BP_{\lambda_1\otimes\nu}$ the (fixed) process $W=(W)_{t\ge0}$ is a
Wiener process as above, independent of $\Phi$.

Let $\bF$ denote the space of all $\R^d$-valued functions
on $\R_+$ equipped with the Kolmogorov product $\sigma$-field.
For $w\in\bF$ and $(t_1,x_1)\in[0,\infty)\times\R^d$
we define $w^{t_1,x_1}\in\bF$ by $w^{t_1,x_1}_t:=w_t+\I\{t\ge t_1\}x_1$.
Clearly the mapping $(w,t_1,x_1)\mapsto w^{t_1,x_1}$ is measurable.
Moreover, if $w\in\bD$ then also $w^{t_1,x_1}\in\bD$.
For any measurable $f\dvtx \bF\rightarrow\R$, the measurable function
$\Delta_{t_1,x_1}f\dvtx \bF\rightarrow\R$ is defined by
%
\begin{equation}
\label{ad} \Delta_{t_1,x_1}f(w):=f\bigl(w^{t_1,x_1}\bigr)-f(w),\qquad w \in
\bF.
\end{equation}
Similarly as at (\ref{differn}), we can iterate this definition
to obtain, for $(t_1,x_1,\ldots,t_n,x_n)\in([0,\infty)\times\R^d)^n$ a function
$\Delta^{n}_{t_1,x_1,\ldots,t_n,x_n}f\dvtx \bF\rightarrow\R$.
Further, we define $\Delta^0f:=f$.
For $s>0$ and $w\in\bF$ let $w^{(s)}\in\bF$ be defined
by $w^{(s)}(t):=w(t\wedge s)$ and let $\tilde{\mathcal{A}}_s$ denote
the $\sigma$-field generated by the mapping $w\mapsto w^{(s)}$.
An $\tilde{\mathcal{A}}_s$-measurable function
$f\dvtx \bF\rightarrow\R$ has the property
that $\Delta_{t,x} f\equiv0$ whenever $t>s$.
Define $\mathcal{A}_s:=\tilde{\mathcal{A}}_s\cap\bD$.

In the next theorem, we consider three L\'evy measures $\nu,\nu',\nu^*$.
We assume that $\nu$ and $\nu'$ are absolutely continuous
with respect to $\nu^*$ with densities $g_\nu$ and $g_{\nu'}$,
respectively, that
satisfy
%
\begin{eqnarray}
\label{222} \int(1-g_\nu)^2 \,\mrmd \nu^*+
\int(1-g_{\nu'})^2 \,\mrmd \nu^*&<&\infty,
\\
\label{235}
\int_{|x|\le1}|x|\bigl|1-g_\nu(x)\bigr|\nu^*(\mrmd x)+
\int_{|x|\le1}|x|\bigl|1-g_{\nu'}(x)\bigr|\nu^*(\mrmd x)&<&\infty.
\end{eqnarray}
We also consider $b,b',b^*\in\R^d$ such that
%
\begin{equation}
\label{223} b=b^*+\int_{|x|\le1}x\bigl(g_\nu(x)-1
\bigr)\nu^*(\mrmd x),\qquad b'=b^*+\int_{|x|\le1}x
\bigl(g_{\nu'}(x)-1\bigr)\nu^*(\mrmd x).
\end{equation}
In the following theorem and also later,
we abuse our notation by interpreting for
a function $g\dvtx \R^d\rightarrow\R$ and $n\in\N$, $g^{\otimes n}$
as a function on $([0,\infty)\times\R^d)^n$.
%
\begin{theorem}\label{tvarlevy}
Assume that (\ref{222}), (\ref{235}) and (\ref{223}) hold.
Let $f\dvtx \bD\rightarrow\R$ be $\mathcal{A}_{t_0}$-measurable
for some $t_0>0$ and assume that $\BE_{b^*,\nu^*}f(X)^2<\infty$.
Then
%
\begin{equation}
\label{mainlevy} \sum^\infty_{n=0}
\frac{1}{n!}\int\bigl|\BE_{b,\nu} \Delta^n f(X)\bigr|
|g_{\nu'}-g_{\nu}|^{\otimes n} \,\mrmd \bigl(\lambda_1
\otimes\nu^*\bigr)^n<\infty,
\end{equation}
where $\BE_{b,\nu} \Delta^n f(X)$ denotes the function
$(t_1,x_1,\ldots,t_n,x_n)\mapsto\BE_{b,\nu} \Delta^n_{t_1,x_1,\ldots
,t_n,x_n}f(X)$.
Furthermore,
%
\begin{equation}
\label{mainlevy1} \BE_{b',\nu'} f(X)=\BE_{b,\nu} f(X) +\sum
^\infty_{n=1}\frac{1}{n!}\int\bigl(
\BE_{b,\nu} \Delta^n f(X)\bigr) (g_{\nu'}-g_{\nu})^{\otimes n}\,\mrmd
\bigl(\lambda_1\otimes\nu^*\bigr)^n.
\end{equation}
\end{theorem}
\begin{pf} Let $\BX:=[0,\infty)\times\R^d$ and define $\bN$ as before.
Let $\bN_0$ be the measurable set of all $\varphi\in\bN$
such that $\varphi([0,s]\times\{x\dvtx 1/n\le|x|\le n\})<\infty$ for all
$s>0$ and $n\in\N$. Since $\nu$ is a L\'evy measure we have
$\BP_{\lambda_1\otimes\nu}(\Phi\in\bN_0)=1$.
For $\varphi\in\bN_0$ and $n\in\N$, we define
$T^n(\varphi)\in\bF$ by the pathwise integrals
\[
T^n(\varphi)_t:=bt+\int_{1/n\le|x|\le1}\int
^t_0x\bigl(\varphi (\mrmd s,\mrmd x)-\mrmd s\nu(\mrmd x)\bigr) +\int
_{n\ge|x|> 1}\int^t_0x
\varphi(\mrmd s,\mrmd x).
\]
Define $T_{b,\nu}(\varphi)\in\bF$ by
\[
T_{b,\nu}(\varphi)_t:=\liminf_{n\to\infty}
T^n(\varphi),\qquad t\ge0,
\]
whenever this is finite, and by $T_{b,\nu}(\varphi)_t:=0$, otherwise.
For $\varphi\notin\bN_0$ we let $T_{b,\nu}(\varphi)\equiv0$.
Then $T_{b,\nu}$ is a measurable mapping from $\bN$ to $\bF$.
It is a basic property of Poisson and L\'evy processes
(\cite{Kallenberg}, Chapter 15)
that $T^n(\Phi)_t$ converges in $\BP_{\lambda_1\otimes\nu}$-probability
and that
%
\begin{equation}
\label{943} \BP_{\lambda_1\otimes\nu}\bigl(W+T_{b,\nu}(\Phi)\in\cdot\bigr)=
\BP_{b,\nu}(X\in\cdot) \qquad\mbox{on $\bF$},
\end{equation}
where here and later we interpret $X$ also as a random element in $\bF$.
Assumptions (\ref{235}) and (\ref{223}) imply that
$T_{b,\nu}=T_{b',\nu'}=T_{b^*,\nu^*}=:T$, so that the following holds
on $\bF$:
%
\begin{equation}
\label{944} \BP_{\lambda_1\otimes\nu'}\bigl(W+T(\Phi)\in\cdot\bigr)=
\BP_{b',\nu
'}(X\in\cdot),\qquad \BP_{\lambda_1\otimes\nu^*}\bigl(W+T(\Phi)\in\cdot\bigr)=
\BP_{b^*,\nu
^*}(X\in\cdot).
\end{equation}
Let $\lambda_1^{t_0}$ be the restriction of $\lambda_1$ to
the interval $[0,t_0]$. Let $\tilde f\dvtx \bF\rightarrow\R$
be an $\tilde{\mathcal{A}}_{t_0}$-measurable function
satisfying $\BE_{b^*,\nu^*}\tilde f(X)^2<\infty$.
We apply Theorem \ref{tvar3} and Remark \ref{rrand}
with $(\lambda,\nu,\rho)$ replaced with
$(\lambda_1^{t_0}\otimes\nu,\lambda_1^{t_0}\otimes\nu',\lambda_1^{t_0}\otimes\nu^*)$,
with $\eta=W$ and
 with the function $(w,\varphi)\mapsto\tilde f(w+T(\varphi))$.
Assumption (\ref{22a}) is implied by (\ref{222}), while
$\BE_{\lambda_1^{t_0}\otimes\nu^*}(W+\tilde f(T(\Phi)))^2<\infty$
follows from (\ref{944})
and assumption on $\tilde f$. (By $\tilde{\mathcal
{A}}_{t_0}$-measurability of
$\tilde f$ we have
$\tilde f(T(\varphi))=\tilde f(T(\varphi_{t_0}))$ for any $\varphi
\in\bN_0$,
where $\varphi_{t_0}$ is the restriction of $\varphi$ to
$[0,t_0]\times\R^d$.) Using that for $\varphi\in\bN_0$,
\[
D^n_{(t_1,x_1),\ldots,(t_n,x_n)}(\tilde f\circ T) (\varphi) =\bigl(
\Delta^n_{t_1,x_1,\ldots,t_n,x_n}\tilde{f}\bigr)
\bigl(T(\varphi)\bigr),\quad
(t_1,x_1,\ldots,t_n,x_n)\in
\bigl([0,\infty)\times\R^d\bigr)^n,
\]
we obtain (\ref{mainlevy}) and (\ref{mainlevy1}) with
$\tilde f$ instead of $f$.

To conclude the proof, we need a
$\tilde{\mathcal{A}}_{t_0}$-measurable function
$\tilde f\dvtx \bF\rightarrow\R$ such that $f=\tilde f$ on~$\bD$.
Such a function trivially exists if $f(w)=g(w_{t_1},\ldots,w_{t_n})$,
where $0\le t_1\le\cdots\le t_n\le t_0$ and $g\dvtx \R^n\rightarrow\R$
is Borel-measurable. Therefore, the existence follows by a
monotone class argument.
\end{pf}
%
\begin{remark}\label{r237}
In the above proof, we cannot apply Theorem \ref{tvar3}
with $(\lambda,\nu,\rho)$ replaced with
$(\lambda_1\otimes\nu,\lambda_1\otimes\nu',\lambda_1\otimes\nu^*)$.
For instance, the first integral in (\ref{22a}) would diverge as soon
as $\nu\ne\nu^*$. Therefore, we have assumed the function $f$
to depend only on the restriction of $X$ to a finite time interval.
\end{remark}
%
\begin{remark}\label{r39}
If
%
\begin{equation}
\label{ord} \int\bigl(|x|\wedge1\bigr) \nu(\mrmd x)<\infty,
\end{equation}
it is common, to rewrite (\ref{LCh}) as
%
\begin{equation}
\label{LCh2} X_t=at+W_t+\int_{\R^d}
\int^t_0x \Phi(\mrmd s,\mrmd x),\qquad t\ge0,
\end{equation}
where $a:=b-\int_{|x|\le1} x \nu(\mrmd x)$. If all three measures
$\nu,\nu',\nu^*$ satisfy (\ref{ord}), then we might replace
$(b,b',b^*)$ by $(a,a',a^*)$ (with $a'$ and $a^*$ defined similarly as $a$)
and simplify (\ref{223}) to $a=a'=a^*$.
\end{remark}
%
\begin{remark}\label{r34}
By \cite{JaSh87}, Theorem IV.4.39, the finiteness of the first
integrals in (\ref{222}) and (\ref{235}) together with the first
identity in (\ref{223})
imply that $\BP_{b,\nu}(X^{(t)}\in\cdot)$ is, for every $t\ge0$,
absolutely continuous
with respect to $\BP_{b^*,\nu^*}(X^{(t)}\in\cdot)$. (Recall that
$X^{(t)}_s:=X_{t\wedge s}$.)
In fact, this conclusion remains true under the weaker assumption
$\int(1-\sqrt{g_\nu})^2\,\mrmd \nu^*<\infty$. We do not know whether the
assumption (\ref{222}) in Theorem \ref{tvarlevy} can be weakend to
$\int(1-\sqrt{g_\nu})^2\,\mrmd \nu^*+\int(1-\sqrt{g_{\nu'}})^2\,\mrmd \nu^*<\infty$,
see also Remark \ref{ralt}.
\end{remark}

Our next theorem is the L\'evy version of Theorem \ref{tseries}.
We consider a L\'evy measure $\nu$ with
density $g_\nu$ with respect to some other L\'evy measure $\nu^*$ and
a measurable
function $g\dvtx \R^d\rightarrow\R$ such that
%
\begin{equation}
\label{567} \int(1-g_\nu)^2\,\mrmd \nu^*<\infty,\qquad \int
g^2\,\mrmd \nu^*<\infty,\qquad \int\bigl(\|x\|\wedge1\bigr)\bigl|g(x)\bigr|\nu^*(\mrmd x)<\infty.
\end{equation}

\begin{theorem}\label{tseriesL} Assume that (\ref{567}) holds
and the first integral in (\ref{235}) is finite.
Let $b$ and $b^*$ satisfy the first identity in (\ref{223}).
Let $\theta_0\in\R$ and assume that $I\subset\R$ is an interval
with non-empty interior such that $\theta_0\in I$ and
$g_\theta:=g_\nu+(\theta-\theta_0)g\ge0$ $\nu^*$-a.e. for
$\theta\in I$.
For $\theta\in I$ let
%
\begin{equation}
\label{7123} b_\theta:=b+(\theta-\theta_0)\int
_{|x|\le1}xg(x)\nu^*(\mrmd x)
\end{equation}
and let $\nu_\theta$ denote the measure
with density $g_\theta$ with respect to $\nu^*$.
Let $f\dvtx \bD\rightarrow\R$ be $\mathcal{A}_{t_0}$-measurable
for some $t_0>0$ and assume that $\BE_{b^*,\nu^*}f(X)^2<\infty$.
Then
%
\begin{eqnarray}
\label{pslevy} \BE_{b_\theta,\nu_\theta} f(X)&=&\BE_{b,\nu} f(X)\nonumber\\[-8pt]\\[-8pt]
&&{} +\sum
^\infty_{n=1}\frac{(\theta-\theta_0)^n}{n!}\int\bigl(
\BE_{b,\nu
} \Delta^n f(X)\bigr) g^{\otimes n}\,\mrmd \bigl(
\lambda_1\otimes\nu^*\bigr)^n, \qquad\theta\in I,\nonumber
\end{eqnarray}
where the series converges absolutely. In particular,
%
\begin{equation}
\label{deriv77} \frac{\mrmd}{\mrmd \theta}\BE_{b_\theta,\nu_\theta} f(X) \bigg|_{\theta
=\theta_0} =
\iint\bigl(\BE_{b,\nu} \Delta_{t,x}f(X)\bigr)g(x)\,\mrmd t\nu^*(\mrmd x).
\end{equation}
\end{theorem}
\begin{pf} Noting that
%
\begin{equation}
\label{r512} b_\theta=b^*+\int_{|x|\le1}x
\bigl(g_\theta(x)-1\bigr)\nu^*(\mrmd x),
\end{equation}
and using the mapping $T$ defined in the proof
of Theorem \ref{tvarlevy},
the result follows from Theorem~\ref{tseries}
and Remark \ref{rrand}.
\end{pf}
%
\begin{remark}\label{rgatlev}
Consider $\nu$ and $\nu^*$ such that the first integrals
in (\ref{222}), respectively, in (\ref{235}) are finite.
Let $b$ and $b^*$ satisfy the first identity in (\ref{223}).
Let $f\dvtx \bD\rightarrow\R$ be a measurable function
such that $\BE_{b^*,\nu^*}f(X)^2<\infty$.
By Theorem \ref{tseriesL}
%
\begin{equation}
\label{Gnu} G_{b,\nu,f}(g):=\iint\bigl(\BE_{b,\nu}
\Delta_{t,x}f(X)\bigr)g(x)\,\mrmd t\nu^*(\mrmd x)
\end{equation}
can be interpreted as the G\^ateaux derivative of
the mapping $\nu'\mapsto\BE_{\nu',b'}f(X)$
at $\nu$ in the direction~$g$, where
$b'$ is determined by $b$, and $\nu'$ and the function $g$ satisfies
the second and third equality in (\ref{567})
as well as $g_\lambda+\theta g\ge0$ $\nu^*$-a.e. for
all $\theta$ in some open neighborhood of $0$.
Proposition \ref{rfrech} on Fr\'echet derivatives can be adapted in a
similar way.
Details are left to the reader.
\end{remark}

The next result deals with non-linear perturbations
and is a consequence of Theorem~\ref{tderiv}.
%
\begin{theorem}\label{tderivL} Assume that (\ref{567}) holds.
Let $\theta_0\in\R$ and assume that $I\subset\R$ is an interval
with non-empty interior such that $\theta_0\in I$. For any $\theta\in
I$ let
$R_\theta\dvtx \R^d\rightarrow[0,\infty)$ be a measurable function
such that the following assumptions are satisfied:
\begin{longlist}[(iii)]
\item[(i)] For all $\theta\in I$,
$g_\nu+(\theta-\theta_0)(g+R_\theta)\ge0$ $\nu^*$-a.e.
\item[(ii)]
$\int(|x|\wedge1) |R_\theta(x)|\nu^*(\mrmd x)<\infty$.
\item[(iii)] $\lim_{\theta\to\theta_0}R_\theta=0$ $\nu^*$-a.e.
\item[(iv)] There is a measurable function $R\dvtx \R^d\rightarrow
[0,\infty)$
such that $|R_\theta|\le R$ $\nu^*$-a.e. and $\int R(x)^2\nu^*(\mrmd x)<\infty$.
\end{longlist}
For $\theta\in I$, let $\nu_\theta$ denote the measure
with density $g_\nu+(\theta-\theta_0)(g+R_\theta)$ with respect to
$\nu^*$.
Let $b,b^*\in\R$ satisfy the first identity in (\ref{223}) and define
%
\begin{equation}
\label{bb} b_\theta:=b+(\theta-\theta_0)\int
_{|x|\le1}x\bigl(g(x)+R_\theta(x)\bigr)\nu^*(\mrmd x).
\end{equation}
Let $f\dvtx \bD\rightarrow\R$ be $\mathcal{A}_{t_0}$-measurable
for some $t_0>0$ and such that
$\BE_{b^*,\nu^*}f(X)^2<\infty$. Then (\ref{deriv77}) holds.
\end{theorem}
%
\begin{remark}\label{rderivL}
Assume that we can take $\nu^*=\nu$ in Theorem \ref{tderivL} (yielding
that $\nu_\theta\ll\nu$). By Theorem \ref{t57}, assumptions (iii) and
(iv) can then be replaced with $\lim_{\theta\to\theta_0}\int R_\theta^2
\,\mrmd \nu=0$.
\end{remark}

We finish this section with some examples.
%
\begin{example}\label{exst}
Let $\alpha\in(0,2)$ and let $\BQ$ be a finite measure on the unit
sphere $\bS:=\{x\in\R^d\dvtx |x|\le1\}$. Then
\[
\nu:=\int_{\bS}\int^\infty_0
\I\{ru\in\cdot\}r^{-\alpha
-1}\,\mrmd r\BQ(\mrmd u)
\]
is the L\'evy measure of an \textit{$\alpha$-stable} L\'evy process,
see, for example, \cite{Bertoin96}. Consider the
L\'evy measure
\[
\mu:=\int_{\bS}\int^1_0 \I
\{ru\in\cdot\}r^{-\alpha'-1}\,\mrmd r\BQ'(\mrmd u),
\]
where $0<\alpha'<\alpha/2$ and $\BQ'$ is a finite measure
on $\bS$. Assume that $\BQ'\ll\BQ$ with a density that is
square-integrable with respect to $\BQ$.
It is not difficult to check
that the density $g:=\mrmd \mu/\mrmd \nu$
satisfies the assumptions of Theorem \ref{tseriesL} with $\nu^*=\nu$,
$I=[0,\infty)$ and $\theta_0=0$.
\end{example}
%
\begin{example}\label{exgamma1}
Let $d=1$, $\alpha\in(0,2)$ and $\nu(\mrmd x):=\I\{x\ne0\} x^{-\alpha-1}\,\mrmd x$
be the L\'evy measure of a (symmetric) $\alpha$-stable process. It is
again easy to check that, for $\beta>0$, the density $g$ of the measure
$\mu_{\beta}(\mrmd x):=\I\{x>0\} x^{-1}\RMe^{-\beta x}\,\mrmd x$ with respect to $\nu$
satisfies the assumptions of Theorem~\ref{tseriesL} with $\nu^*=\nu$,
$I=[0,\infty)$ and $\theta_0=0$. For $\theta\ge0$, the measure
$\theta\mu$ is the L\'evy measure of a \textit{Gamma process} with shape
parameter $\theta$ and scale parameter $\beta$, see, for example,
\cite{Bertoin96}. (Under $\BP_{0,\theta\mu_\beta}$ and for $W\equiv0$,
the random variable $X_t$ has a Gamma distribution with shape parameter
$\theta t$ and scale parameter $\beta$.)
\end{example}
%
\begin{example}\label{exgamma2}
We let $\nu$ and $\mu_{\beta}$ be as in Example \ref{exgamma1}. This
time we are interested in derivatives with respect to the scale
parameter $\beta$. We fix $\beta_0>0$ and $\theta\ge0$. Our aim is to
apply Theorem \ref{tderivL} with $I=(\beta_0/2,3\beta_0/2)$,
$\theta_0=\beta_0$, $\nu_\beta:= \nu+\theta\mu_{\beta}$, and
$\nu^*:=\nu$. The measure $\nu$ in Theorem \ref{tderivL} is being
replaced with $\nu_{\beta_0}$. We have noted in Example \ref{exgamma1}
that $\int(1-g_{\beta_0})^2\,\mrmd \nu<\infty$, where $g_{\beta_0}$ is the
Radon--Nikodym derivative of $\nu_{\beta_0}$ with respect to $\nu$.
Next, we note that
%
\begin{eqnarray}
\label{315} \nu_\beta(\mrmd x) &=& \nu_{\beta_0}(\mrmd x)+\I\{x>0\}\theta
\bigl(\RMe^{-(\beta-\beta
_0)x}-1\bigr)x^{\alpha}\RMe^{-\beta_0 x}\nu(\mrmd x)
\nonumber\\[-8pt]\\[-8pt]
&=& \bigl[g_{\beta_0}(x)+(\beta-\beta_0) \bigl(g(x)+R_\beta(x)
\bigr)\bigr]\nu(\mrmd x),\nonumber
\end{eqnarray}
where
\[
g(x):=-\I\{x>0\}\theta x^{\alpha+1}\RMe^{-\beta_0 x},\qquad R_\beta(x):=-
\I\{x>0\}\theta x^{\alpha}\RMe^{-\beta_0 x} \sum
^\infty_{n=2}\frac{(\beta_0-\beta)^{n-1}}{n!}x^n.
\]
Clearly, $g(x)^2$ and $(|x|\wedge1)|g(x)|$ are integrable
with respect to $\nu$.
We need to check the assumptions (i)--(iv) of Theorem \ref{tderivL}.
Assumption (i) follows from (\ref{315}) while (iii) is obvious.
Further, we have for $\beta\in I$ that $|R_\beta|\le R$, where
\[
R(x):=\I\{x>0\}\theta x^{\alpha}\RMe^{-\beta_0 x}\sum
^\infty_{n=0} \frac{(\beta_0/2)^{n-1}}{n!}x^n =\I\{x>0
\}\theta c^{-1} x^{\alpha}\RMe^{-c x},
\]
where $c:=\beta_0/2$. Since $\int(x\wedge1) R(x)\nu(\mrmd x)<\infty$
and $\int R(x)^2\nu(\mrmd x)<\infty$ we obtain (ii) and (iv).
Let $b\in\R$. In view of (\ref{r512}), we define
\[
b_\beta:=b+\int^1_0x
\bigl(g_\beta(x)-1\bigr)\nu(\mrmd x) =b+\theta\int^1_0
\RMe^{-\beta x} \,\mrmd x =b+\frac{\theta}{\beta} \bigl(1-\RMe^{-\beta} \bigr).
\]
Under the assumption $\BE_{b,\nu}f(X)^2<\infty$,
we then obtain that
%
\begin{equation}
\label{gammaderivscale} \frac{\mrmd}{\mrmd \beta}
\BE_{b_{\beta},\nu+\theta\mu_{\beta}}f(X) \bigg|_{\beta=\beta_0}
=-\theta\iint\bigl(\BE_{b_{\beta_0},\nu+\theta\mu_{\beta_0}}\Delta_{t,x}f(X)
\bigr)\RMe^{-\beta_0 x}\,\mrmd t \,\mrmd x.
\end{equation}
\end{example}
%
\begin{remark}\label{r42}
It is a common feature of Examples \ref{exst} and \ref{exgamma1} that
the perturbation $\nu'-\nu$ is infinite. Theorem \ref{tvar} would not
be enough to treat these cases. In Example \ref{exgamma2} however,
$\nu_{\beta_0}-\nu_\beta$ is finite so that one might use Theorem
\ref{tvar} in the case $\beta<\beta_0$.
\end{remark}

Finally in this section we assume $d=1$ and apply our results to the
\textit{running supremum}
\[
S_t:=\sup\{X_s\dvtx 0\le s\le t\},\qquad t\ge0,
\]
of $X$. We fix $t_0>0$ and define
\[
Z_t:=\sup\{X_s\dvtx t\le s\le t_0\},\qquad
Y_t:=S_t-Z_t,\qquad 0\le t\le t_0.
\]

\begin{proposition}\label{p6} Let $\nu,\nu^*,g,g_\nu,b,b^*,I,\theta_0$ be as in
Theorem \ref{tseriesL} and define $b_\theta,g_\theta$ as in that
theorem. Assume moreover that
%
\begin{equation}
\label{mom5} \int_{x>1}x^2\nu^*(\mrmd x)<\infty.
\end{equation}
Then $\theta\mapsto\BE_{b_{\theta},\nu_{\theta}} S_{t_0}$
is analytic on $I$. Moreover,
%
\begin{equation}
\label{supderiv} \frac{\mrmd}{\mrmd \theta} \BE_{b_{\theta},\nu_{\theta}}S_{t_0}
\bigg|_{\theta=\theta_0} =\iint\bigl((x-y)^+-y^-\bigr)g(x)\BQ(\mrmd y) \nu^*(\mrmd x),
\end{equation}
where, for $x\in\R$, $x^+:=\max\{x,0\}$, $x^-:=-\min\{x,0\}$, and
%
\begin{equation}
\label{aa} \BQ:=\int^{t_0}_0\BP_{b,\nu}(Y_t
\in\cdot)\,\mrmd t.
\end{equation}
\end{proposition}
\begin{pf} We define a measurable function $f\dvtx \bD\rightarrow\R$ by
\[
f(w):=\sup\{w_s\dvtx 0\le s\le t_0\},\qquad w\in\bD.
\]
It follows from the L\'evy--Khintchine representation (\ref{LCh}),
Doob's inequality and moment properties of Poisson integrals
that (\ref{mom5}) is sufficient (and actually also necessary)
for \mbox{$\BE_{b^*,\nu^*}f(X)^2<\infty$}. (This argument is quite standard.)
Hence, we can apply Theorem \ref{tseriesL}.

It remains to compute the right-hand side of (\ref{deriv77}).
Let $t\in(0,t_0]$. For $x> 0$ we have
\[
f\bigl(X^{t,x}\bigr)= \cases{ S_{t_0}+x, &\quad if
$S_{t-}\le Z_t$,
\cr
Z_{t}+x, &\quad if
$Z_t<S_{t-}\le Z_t+x$,
\cr
S_{t_0}, &\quad
if $Z_t+x<S_{t-}$, }
\]
so that
\[
f\bigl(X^{t,x}\bigr)-f(X)=\I\{Y_t\le0\}x+\I
\{0<Y_t\le x\}(x-Y_t) =(x-Y_t)^+-(Y_t)^-,
\]
provided that $S_{t-}=S_t$. Note that the latter equality holds for
$\lambda_1$-a.e. $t>0$. Similarly we obtain for $x<0$ that
\[
f\bigl(X^{t,x}\bigr)-f(X)=\I\{x<Y_t\le0\}Y_t+\I
\{Y_t\le x\}x=(x-Y_t)^+-(Y_t)^-,
\]
whenever $S_{t-}=S_t$.
Hence, (\ref{supderiv}) follows from (\ref{deriv77}).
Note that the integrability required for (\ref{supderiv}) is part
of the assertion of Theorem \ref{tseriesL}. But it does also follow
more directly from $|f(X^{t,x})-f(X)|\le2 |x|$ and the fact that
$\int|g(x)||x|\nu^*(\mrmd x)$ is finite by assumption (\ref{567}) on
$g$, (\ref{mom5}), and the Cauchy--Schwarz inequality.
\end{pf}
%
\begin{remark}\label{r714} We consider the situation of
Proposition \ref{p6} but do not assume (\ref{mom5}).
For $u\in\R$ wen can then apply Theorem \ref{tseriesL} to
the real and the imaginary part of the
complex-valued and bounded function $f(X):=\RMe^{iuS_{t_0}}$.
This shows that $\theta\mapsto\BE_{b_{\theta},\nu_{\theta}} f(X)$
is analytic. The derivative at $\theta_0$ can be expressed
in terms of the measure
%
\begin{equation}
\label{bbp} \BQ':=\int^{t_0}_0
\BP_{b,\nu}\bigl((S_t,Z_t)\in\cdot\bigr)\,\mrmd t
\end{equation}
that contains more information than the measure (\ref{aa}).
The same remark applies to the bounded function
$f(X):=\I\{S_{t_0}>u\}$. The details are left to the reader.
\end{remark}

For a general L\'evy process the distribution of $S_t$ is not
known. The measures (\ref{aa}) and (\ref{bbp}) are not known
either. This hints at the fact that perturbation analysis
cannot help in finding explicit distributions.
What equation (\ref{supderiv}) does, however, is to identify the
G\^ateaux derivative of $\BE_{b,\nu} S_{t_0}$ in the direction $g$,
see Remark \ref{rgatlev}. The measure (\ref{aa}), controlling all
these derivatives, is completely determined by the distribution
of the process $(X_t)_{t\le t_0}$ under $\BP_{b,\nu}$.
We do not make any attempt to review the vast literature
on the running supremum of L\'evy processes but just refer to
\cite{Chaumont12} for some recent progress.

\section*{Acknowledgements}

The paper benefited from several very useful comments and proposals of
two referees and the Associated Editor.



\printhistory

\end{document}